\documentclass[11pt]{article}
\usepackage{latexsym,amssymb,amsmath,a4wide,theorem}
\usepackage[isolatin]{inputenc}
\usepackage[T1]{fontenc}
\usepackage{times}
\usepackage{cite}
\usepackage{amsfonts}
\usepackage{hyperref}

\def\sqr#1#2{{\vbox{\hrule height.#2pt
     \hbox{\vrule width.#2pt height#1pt \kern#1pt
           \vrule width.#2pt}
     \hrule height.#2pt}}}

\newtheorem{theorem}{Theorem}[section]
\newtheorem{lemma}[theorem]{Lemma}
\newtheorem{remark}[theorem]{Remark}

\title{\Large\bf Concentration of positive ground state solutions for critical Kirchhoff equation with competing potentials}
\author{
{Yongpeng Chen$^{1}$},
{Zhipeng Yang$^{2}$}\thanks{zhipeng.yang@mathematik.uni-goettingen.de},\\
\small School of Science, Guangxi University of Science and Technology, Liuzhou 545006. P.R.China.$^{1}$\\
\small Mathematical Institute, Georg-August-University of G\"ottingen, G\"ottingen 37073, Germany.$^{2}$\\
}

\date{}

\begin{document}
\maketitle
\begin{abstract}
In this paper, we consider the following singularly perturbed Kirchhoff equation
\begin{equation*}
-(\varepsilon^2a+\varepsilon b\int_{\mathbb{R}^3}|\nabla u|^2dx)\Delta u+V(x)u=P(x)|u|^{p-2}u+Q(x)|u|^4u,\quad x\in\mathbb{R}^3,
\end{equation*}
where $\varepsilon>0$ is a small parameter, $a, b > 0$ are constants, $p\in(4,6)$ and $V, P, Q$ are potential functions satisfying some competing conditions. We prove the existence of a positive ground state solution by using variational methods, and we determine a concrete set related to the potentials $V,P$ and $Q$ as the concentration position of these ground state solutions as $\varepsilon\to0$.
\end{abstract}


{\bf AMS:} Subject Classification: 35Q55,  35J655
\numberwithin{equation}{section}
\section{Introduction and Main Results}\
\par
In this paper, we investigate the existence and concentration behavior of positive ground
solutions to the following Kirchhoff type equation with critical exponent
\begin{equation}\label{maineq1}
-(\varepsilon^2a+\varepsilon b\int_{\mathbb{R}^3}|\nabla u|^2dx)\Delta u+V(x)u=P(x)|u|^{p-2}u+Q(x)|u|^4u,\quad x\in\mathbb{R}^3,
\end{equation}
where $\varepsilon>0$ is a small parameter, $a, b > 0$ are constants, $p\in(4,6)$.  This problem motivated by some works related to the following Kirchhoff equation
\begin{equation}\label{eq1.2}
\begin{cases}
-M\bigg(\displaystyle{\int_{\Omega}}|\nabla u|^2dx\bigg)\Delta u=f(x,u)&\quad \text{in} \quad\Omega,\\
u=0&\quad\text{on}\quad\partial\Omega,
\end{cases}
\end{equation}
where $M(t)=a+bt$ ($a,b>0$) for all $t\geq0$. Such a problem is often referred to as being nonlocal because of the presence of the term $M(\int_{\Omega}|\nabla u|^2dx)$ which implies that the equation \eqref{eq1.2} is no longer a pointwise identity. This phenomenon leads to some mathematical difficulties, which make the study of such a class of problems particularly interesting,  see for example \cite{Arosio-Panizzi1996TAMS} and \cite{DpAncona-Spagnolo1992IM} for more information about \eqref{eq1.2}. Recall that \eqref{maineq1} is called degenerate when $a=0$ and $b>0$, and a nondegenerate one when $a>0$ and $b>0$ (see e.g. \cite{DpAncona-Spagnolo1992IM}, \cite{Ono1997JDE}).
\par
On one hand, \eqref{eq1.2} is related to the stationary analogue of the Kirchhoff equation
\begin{equation}\label{eq1.3}
\begin{cases}
u_{tt}-M\bigg(\displaystyle{\int_{\Omega}}|\nabla u|^2dx\bigg)\Delta u=f(x,u)&\quad\text{in} \quad\Omega,\\
u=0&\quad\text{on}\quad\partial\Omega.
\end{cases}
\end{equation}
It was proposed by Kirchhoff in \cite{Kirchhoff1883} as a generalization of the well-known D'Alembert wave equation
$$\rho\frac{\partial^2u}{\partial t^2}-\bigg(\frac{p_0}{\lambda}+\frac{E}{2L}\int_0^L\bigg|\frac{\partial u}{\partial x}\bigg|^2dx\bigg)\frac{\partial^2u}{\partial x^2}=f(x,u)$$
for free vibrations of elastic strings. It seems that the first result concerning the global sovability was proved by Bernstein \cite{Bernstein1940BASUS}. This result was generalized to the case $n\geq 1$ by  Pohozaev in \cite{Pohozaev1975MS}. In \cite{Lions1978NHMS}, Lions proposed an an abstract framework to the problem. Since then, \eqref{eq1.3} received much more attention. We have to point out that nonlocal problems also appear in other fields as biological systems, where $u$ describes a process which depends on the average of itself (for example, population density). See, for example, \cite{Arosio-Panizzi1996TAMS} and the references therein. To the best of our knowledge, the variational methods were first involved in \cite{Alves-Correa-Ma2005CMP} and \cite{Ma-Rivera2003AML}. After that, there have been many works about the existence of nontrivial solutions to \eqref{eq1.2} by using different variational techniques, see e.g. \cite{Chen-Kuo-Wu2011JDE,Deng-Peng-Shuai2015JFA,He-Zou2009NA,Mao-Zhang2009NL,Perera-Zhang2006JDE,Shuai2015JDE,Wu2011NA,Zhang-Perera2006JMAA} and the references therein.
\par
On the other hand, \eqref{maineq1} can come back to the following equation
\begin{equation}\label{eq1.4}
\begin{cases}
-\bigg(\varepsilon^2a+ \varepsilon b\displaystyle{\int_{\mathbb{R}^3}}|\nabla u|^2dx\bigg)\Delta u+ V(x)u=f(u)&\quad x\in\mathbb{R}^3,\\
u\in H^1(\mathbb{R}^3)&\quad x\in\mathbb{R}^3.
\end{cases}
\end{equation}
The existence and multiplicity of solutions to \eqref{eq1.4} with $\varepsilon=1$ were studied in some recent works. Li and Ye \cite{Li-Ye2014JDE} obtained the existence of a positive ground state of \eqref{eq1.4} with $f(u)=|u|^{p-1}u$ for $2<p<5$.
In \cite{Deng-Peng-Shuai2015JFA}, Deng, Peng and Shuai studied the existence and asymptotical behavior of nodal solutions of \eqref{eq1.4} with $V$ and $f$ is radially symmetric in $x$ as $b\to 0^+$.
Very recently, Li et.al proved that the positive ground state solution of \eqref{eq1.4}
with $V\equiv1$ and $f(u)=|u|^{p-1}u$ ($1<p<5$) is unique and nondegenerate (see \cite{MR4021897}).
\par
For the concentration behavior of solutions as $\varepsilon\to 0^+$, He and Zou \cite{He-Zou2012JDE} considered the multiplicity and concentration behavior of the positive solutions of \eqref{eq1.4} by using Ljusternik-Schnirelmann theory (see \cite{Willem1996}) and minimax methods, the author obtained the multiplicity of positive solutions, which concentrate on the global minima of $V(x)$ as $\varepsilon\rightarrow 0^+$. A similar result for the critical case $f(u)=\lambda g(u)+|u|^4u$ was obtained separately in \cite{He-Zou2014ADM} and \cite{Wang-Tian-Xu-Zhang2012JDE}, where the subcritical term $g(u)\sim |u|^{p-2}u$ with $4<p<6$.
In \cite{He-Li-Peng2014ANS}, He, Li and Peng constructed a family of positive solutions which
concentrates around a local minimum of $V$ as $\varepsilon\to 0^+$ for a critical problem $f(u)=g(u)+|u|^4u$ with $g(u)\sim |u|^{p-2}u$ ($4<p<6$). For the more delicate case that  $f(u)=\lambda |u|^{p-2}u+|u|^4u$ with $2<p\leq 4$ we refer to He and Li \cite{He-Li2015CVPDE}, where a family of positive solutions which
concentrates around a local minimum of $V$ as $\varepsilon\to 0^+$ were obtained.
\par
In this paper, we are concerned with the existence and concentration behavior of ground state solutions for \eqref{maineq1}.
We note that \eqref{maineq1} involves  critical exponent and three different potentials which make our problem more complicated.
This brings a competition between the potentials $V$, $P$ and $Q$: each one would like to attract ground states to their minimum or maximum points,
respectively. It makes difficulties in determining the concentration position of solutions. This kind of problem can be traced back to \cite{Wang1993CMP}, and \cite{Wang-Zeng1997SIAMJMA} for the semilinear Schr\"{o}dinger equation.
In \cite{Ding-Liu2012JDE}, the authors found new concentration phenomena for Dirac equations with competing potentials and subcritical or critical nonlinearities, respectively. See also \cite{Ding-Liu2013MM,MR3518335,MR3994307} for other related results.
\par
We need some notations to help us to determine the concentration set of solutions. Set
\begin{equation*}
\begin{split}
&0<V_{\min}:=\min_{x\in \mathbb{R}^3} V(x),~~V_{\max}:=\sup_{x\in \mathbb{R}^3}V(x),~~\mathcal{V}:=\{x\in \mathbb{R}^3:V(x)=V_{\min}\},~~V_\infty:=\liminf_{|x|\rightarrow \infty}V(x),\\
&0<P_{\min}:=\inf_{x\in \mathbb{R}^3}P(x),~~P_{\max}:=\max_{x\in \mathbb{R}^3}P(x),~~\mathcal{P}:=\{x\in \mathbb{R}^3:P(x)=P_{\max}\},\ P_\infty:=\limsup_{|x|\rightarrow \infty}P(x),\\
&0<Q_{\min}:=\inf_{x\in \mathbb{R}^3}Q(x),\ Q_{\max}:=\max_{x\in \mathbb{R}^3}Q(x),~~\mathcal{Q}:=\{x\in \mathbb{R}^3:Q(x)=Q_{\max}\},\ Q_\infty:=\limsup_{|x|\rightarrow \infty}Q(x),\\
\end{split}
\end{equation*}
Moreover, we assume that $V,P$ and $Q$ are three locally H\"older continuous and bounded functions satisfying the following conditions:
\begin{itemize}
\item[$(PQ1)$]$\mathcal{P\cap Q}=\{x\in\mathbb{R}^3:  P(x)=P_{max}, Q(x)=Q_{max}\}\neq \varnothing.$
\item[$(PQ2)$]$P_{max}>P_{\infty}$ and there exist $R>0$ and $x^\ast\in \mathcal{P\cap Q}$ such that $V(x^\ast)\leq V(x)$ for all $|x|\geq R.$
\item[$(VQ1)$]$\mathcal{V\cap Q}=\{x\in\mathbb{R}^3:  V(x)=V_{min}, Q(x)=Q_{max}\}\neq \varnothing.$
\item[$(VQ2)$]$V_{\infty}>V_{min}$ and there exist $R>0$ and $x^\ast\in \mathcal{V\cap Q}$ such that $P(x^\ast)\geq P(x)$ for all $|x|\geq R.$
\end{itemize}
\par
Define the following set
\begin{equation*}
\mathcal{A}_V=\{x\in \mathcal{P\cap Q}: V(x)=V(x^\ast)\}\cup\{x\notin \mathcal{P\cap Q}: V(x)<V(x^\ast)\}
\end{equation*}
and
\begin{equation*}
\mathcal{A}_P=\{x\in \mathcal{V\cap Q}: P(x)=P(x^\ast)\}\cup\{x\notin \mathcal{V\cap Q}: P(x)>P(x^\ast)\}.
\end{equation*}
Obviously, under the assumptions $(PQ1)$ and $(PQ2)$, the set $\mathcal{A}_V$ is bounded and  we can assume $V(x^\ast)=\displaystyle\min_{x\in \mathcal{P\cap Q}} V(x).$ Similarly, under the assumptions $(VQ1)$ and $(VQ2)$, the set $\mathcal{A}_P$ is bounded and  we can assume $P(x^\ast)=\displaystyle\max_{x\in \mathcal{V\cap Q}} P(x).$
\par
Now, we can state our main results as follows.
\begin{theorem} Suppose that the potentials $V(x)$, $P(x)$, $Q(x)$ satisfy conditions $(PQ1)$ and $(PQ2)$. Then for any $\varepsilon>0$ small enough,
problem \eqref{maineq1} has at least one positive ground state solution $u_\varepsilon$. Moreover, if $V(x)$, $P(x)$, $Q(x)$ are uniformly continuous on $\mathbb{R}^3$, then
\begin{itemize}
\item[(1)] there exists a maximum point $x_\varepsilon \in \mathbb{R}^3$ of $u_\varepsilon$ such that $\displaystyle\lim_{\varepsilon\rightarrow0}dist(x_\varepsilon,\mathcal{A}_V)=0$, and there exist some constants $c,C>0$ such that
\begin{equation*}
u_\varepsilon(x)\leq Cexp(-\frac{c}{\varepsilon}|x-x_\varepsilon|).
\end{equation*}
\item[(2)]set $\tilde{u}_\varepsilon(x):=u_\varepsilon(\varepsilon x+\tilde{x}_\varepsilon)$, where $\tilde{x}_\varepsilon$ is a maximum point of $\tilde{u}_\varepsilon$. If $x_\varepsilon\rightarrow x_0$ as $\varepsilon\rightarrow 0$, then up to a subsequence, $\tilde{u}_\varepsilon$
    converges in $H^1(\mathbb{R}^3)$ to a positive ground state solution of
\begin{equation*}
-(\varepsilon^2a+\varepsilon b\int_{\mathbb{R}^3}|\nabla u|^2dx)\Delta u+V(x_0)u=P(x_0)|u|^{p-2}u+Q(x_0)|u|^4u,\quad x\in\mathbb{R}^3.
\end{equation*}
In particular if $\mathcal{V}\cap \mathcal{P}\cap \mathcal{Q}\neq \emptyset$, then $\lim\limits_{\varepsilon\rightarrow 0}dist(x_\varepsilon,\mathcal{V}\cap \mathcal{P}\cap \mathcal{Q})=0$, and up to a subsequence, $\tilde{u}_\varepsilon$ converges in $H^1(\mathbb{R}^3)$ to a positive ground state solution of
\begin{equation*}
-(\varepsilon^2a+\varepsilon b\int_{\mathbb{R}^3}|\nabla u|^2dx)\Delta u+V_{min}u=P_{max}|u|^{p-2}u+Q_{max}|u|^4u,\quad x\in\mathbb{R}^3.
\end{equation*}

\end{itemize}
\end{theorem}

\begin{theorem}  Suppose that the potentials $V(x)$, $P(x)$, $Q(x)$ satisfy conditions $(VQ1)$ and $(VQ2)$. Then for any $\varepsilon>0$ small enough,
problem \eqref{maineq1} has at least one positive ground state solution $u_\varepsilon$. Moreover, if $V(x)$, $P(x)$, $Q(x)$ are uniformly continuous on $\mathbb{R}^3$, then
\begin{itemize}
\item[(1)] there exists a maximum point $x_\varepsilon \in \mathbb{R}^3$ of $u_\varepsilon$ such that $\displaystyle\lim_{\varepsilon\rightarrow0}dist(x_\varepsilon,\mathcal{A}_P)=0$, and there exist some constants $c,C>0$ such that
\begin{equation*}
u_\varepsilon(x)\leq Cexp(-\frac{c}{\varepsilon}|x-x_\varepsilon|).
\end{equation*}
\item[(2)]set $\tilde{u}_\varepsilon(x):=u_\varepsilon(\varepsilon x+\tilde{x}_\varepsilon)$, where $\tilde{x}_\varepsilon$ is a maximum point of $\tilde{u}_\varepsilon$. If $x_\varepsilon\rightarrow x_0$ as $\varepsilon\rightarrow 0$, then up to a subsequence, $\tilde{u}_\varepsilon$
    converges in $H^1(\mathbb{R}^3)$ to a positive ground state solution of
\begin{equation*}
-(\varepsilon^2a+\varepsilon b\int_{\mathbb{R}^3}|\nabla u|^2dx)\Delta u+V(x_0)u=P(x_0)|u|^{p-2}u+Q(x_0)|u|^4u,\quad x\in\mathbb{R}^3.
\end{equation*}
In particular if $\mathcal{V}\cap \mathcal{P}\cap \mathcal{Q}\neq \emptyset$, then $\lim\limits_{\varepsilon\rightarrow 0}dist(x_\varepsilon,\mathcal{V}\cap \mathcal{P}\cap \mathcal{Q})=0$, and up to a subsequence, $\tilde{u}_\varepsilon$ converges in $H^1(\mathbb{R}^3)$ to a positive ground state solution of
\begin{equation*}
-(\varepsilon^2a+\varepsilon b\int_{\mathbb{R}^3}|\nabla u|^2dx)\Delta u+V_{min}u=P_{max}|u|^{p-2}u+Q_{max}|u|^4u,\quad x\in\mathbb{R}^3.
\end{equation*}
\end{itemize}
\end{theorem}
\par
It is worth to note that we will overcome some difficulties. The first one is the appearance of the nonlocal $(\int_{\mathbb{R}^3}|\nabla u|^2dx)\Delta u$, one does not know in general
\begin{equation*}
\int_{\mathbb{R}^3}|\nabla u_n|^2dx\int_{\mathbb{R}^3}\nabla u_n\nabla\varphi dx=\int_{\mathbb{R}^3}|\nabla u|^2dx\int_{\mathbb{R}^3}\nabla u\nabla\varphi dx+o_n(1),\quad\forall \varphi\in H^1(\mathbb{R}^3)
\end{equation*}
and
$$\Big(\int_{\mathbb{R}^3}|\nabla u_n|^2dx\Big)^2 -\Big(\int_{\mathbb{R}^3}|\nabla u|^2dx\Big)^2=\Big(\int_{\mathbb{R}^3}|\nabla u_n-\nabla u|^2dx\Big)^2+o_n(1)$$
from $u_n\rightharpoonup u $ in $H^1(\mathbb{R}^3)$.
The second one is that nonlinearity term is critical, then the embedding $H^1(\mathbb{R}^3)\hookrightarrow L^t(\mathbb{R}^3)$ is not compact for any $t\in(2,6)$ so that the standard variational methods can't be applied directly. Thus, some new technical analysis need to be established.
\par
This paper is organized as follows. In the forthcoming section we collect some necessary preliminary Lemmas which will be used later. In section 3, we study the auxiliary problem of \eqref{maineq1}. In section 4, We we are devoted to main results associated with \eqref{maineq1} and some properties as $\varepsilon\to0^+$.
\par
\textbf{Notation.}~In this paper we make use of the following notations.
\begin{itemize}
\item[$\bullet$] For any $R>0$ and for any $x\in\mathbb{R}^3$, $B_{R}(x)$ denotes the ball of radius $R$ centered at $x$.
\item[$\bullet$]  $L^p(\mathbb{R}^3)$, $1\leq p<+\infty$ denotes the Lebesgue space with  the norm $|u|_p=(\int_{\mathbb{R}^3}|u|^pdx)^{\frac{1}{p}}$.
\item[$\bullet$]  $L^{\infty}(\mathbb{R}^3)$ denotes the Lebesgue space with the norm $|u|_{\infty}=ess \sup|f|$.
\item[$\bullet$] The letters $C,C_i$ stand for positive constants (possibly different from line to line).
\item[$\bullet$]  "$\rightarrow$" for the strong convergence and "$\rightharpoonup$" for the weak convergence.
\item[$\bullet$]  $\mu(A)$ denotes the Lebesgue measure of $A\subset\mathbb{R}^3$.
\end{itemize}

\section{Preliminaries}\
\par
Throughout the paper, we consider the Sobolev space $E=H^1(\mathbb{R}^3)$ with the following standard norm
\begin{equation*}
\|u\|=\Big(\int_{\mathbb{R}^3}\big(|\nabla u|^2+u^2\big)dx\Big)^{\frac{1}{2}}
\end{equation*}
and denote the norm of $D^{1,2}(\mathbb{R}^3)$ by
\begin{equation*}
\|u\|_{D^{1,2}}=\Big(\int_{\mathbb{R}^3}\big(|\nabla u|^2\big)dx\Big)^{\frac{1}{2}}.
\end{equation*}
In the following, we denote by $S$ the best Sobolev constant:
\begin{equation*}
S|u|_6^2\leq \int_{\mathbb{R}^3}|\nabla u|^2dx.
\end{equation*}
Making the change of variable $x\mapsto\varepsilon x$ and  $v(x)=u(\varepsilon x)$, problem \eqref{maineq1} reduces to the equation
\begin{equation}\label{transeq1}
-(a+b\int_{\mathbb{R}^3}|\nabla v|^2dx)\Delta v+V(\varepsilon x)v=P(\varepsilon x)|v|^{p-2}v+Q(\varepsilon x)|v|^4v,\quad x\in\mathbb{R}^3.
\end{equation}
Thus, it suffices to study \eqref{transeq1} and the norm
\begin{displaymath}
\|v\|_{\varepsilon}=\big(a\int_{\mathbb{R}^3}|\nabla v|^2dx+\int_{\mathbb{R}^3}V(\varepsilon x)v^2dx\big)^{\frac{1}{2}}
\end{displaymath}
is an equivalent norm on $E$.
\par
The corresponding energy functional
\begin{equation*}
J_\varepsilon(v)=\frac{1}{2}\|v\|^2_\varepsilon+\frac{b}{4}(\int_{\mathbb{R}^3}|\nabla v|^2dx)^2-\frac{1}{p}\int_{\mathbb{R}^3}P(\varepsilon x)|v|^{p}dx-\frac{1}{6}\int_{\mathbb{R}^3}Q(\varepsilon x)|v|^{6}dx.
\end{equation*}
It is easy to check that $J_\varepsilon$ is well defined on $E$ and $J_\varepsilon \in C^1(E,\mathbb{R}).$
\par
Let us define the Nehari manifold \cite{Willem1996} associated with $J_{\varepsilon}$
\begin{equation*}
\mathcal{N}_\varepsilon:=\Big\{u\in E\backslash\{0\}|~ I_\varepsilon(u)=0\Big\},
\end{equation*}
where $I_\varepsilon(u)=\langle J^{\prime}_{\varepsilon}(u),u\rangle$.

\begin{lemma}\label{le:lowboundedinN}There exists $\sigma>0$ which is independent of $\varepsilon$ such that
\begin{equation*}
\|v\|_\varepsilon>\sigma \quad\text{and}\quad J_\varepsilon(v)\geq\frac{p-2}{2p}\sigma^2 \quad\text{for all}\quad v\in  \mathcal{N}_\varepsilon.
\end{equation*}
\end{lemma}

\noindent{\bf Proof:} For any $v\in \mathcal{N}_\varepsilon,$ we have
\begin{equation*}
\begin{aligned}
0&=\|v\|^2_\varepsilon+b(\int_{\mathbb{R}^3}|\nabla v|^2dx)^2-\displaystyle\int_{\mathbb{R}^3}
P(\varepsilon x)|v|^pdx-\int_{\mathbb{R}^3}Q(\varepsilon x)v^6dx.\\
&\geq\|v\|^2_\varepsilon+b(\int_{\mathbb{R}^3}|\nabla v|^2dx)^2-C(\|v\|^p_\varepsilon+\|v\|^6_\varepsilon)
\end{aligned}
\end{equation*}
which implies that there exists $\sigma>0$ such that $\|v\|_\varepsilon>\sigma>0$.
In the above inequality, we have used the boundness of $P(x)$ and $Q(x)$, and the Sobolev embedding Theorem.
\par
On the other hand, we have
\begin{equation*}
\begin{aligned}
J_\varepsilon(v)&=\displaystyle\frac{1}{2}\|v\|^2_\varepsilon+\frac{b}{4}(\int_{\mathbb{R}^3}|\nabla v|^2dx)^2-\frac{1}{p}\int_{\mathbb{R}^3}P(\varepsilon x)|v|^{p}dx-\frac{1}{6}\int_{\mathbb{R}^3}Q(\varepsilon x)|v|^{6}dx,\\
&\geq\displaystyle\frac{1}{2}\|v\|^2_\varepsilon+\frac{b}{4}(\int_{\mathbb{R}^3}|\nabla v|^2dx)^2-\frac{1}{p}\int_{\mathbb{R}^3}P(\varepsilon x)|v|^{p}dx-\frac{1}{p}\int_{\mathbb{R}^3}Q(\varepsilon x)|v|^{6}dx,\\
&=\displaystyle\frac{1}{2}\|v\|^2_\varepsilon+\frac{b}{4}(\int_{\mathbb{R}^3}|\nabla v|^2dx)^2-\frac{1}{p}(\|v\|^2_\varepsilon+b(\int_{\mathbb{R}^3}|\nabla v|^2dx)^2),\\
&\geq\displaystyle(\frac{1}{2}-\frac{1}{p})\|v\|^2_\varepsilon\\
&\geq\displaystyle\frac{p-2}{2p}\sigma^2.
\end{aligned}
\end{equation*}
\hfill{$\Box$}
\begin{remark}
By a direct computation we have $\langle I_\varepsilon^\prime(u),u\rangle<0$, which implies that $I_\varepsilon^\prime(u)\neq0 ,~\forall u\in \mathcal{N}_\varepsilon$. It follows from the Implicity Function Theorem  that $\mathcal{N}_\varepsilon$ is a $C^1$-manifold.
\end{remark}
\par
One can easily check that the functional $J_\varepsilon$ satisfies the mountain-pass geometry, that is the following lemma holds.
\begin{lemma}\label{le:moutain-pass geometry}
$J_\varepsilon$ has the mountain geometry structure.\\
\begin{itemize}
\item[(1)] There exist $a_0,r_0>0$ independent of $\varepsilon$, such that $J_\varepsilon(v)\geq a_0$, for all $v\in E$ with $\|v\|=r_0.$\\
\item[(2)] For any $v\in E\setminus\{0\}$ , $\lim\limits_{t\to \infty}J_\varepsilon(tv)=-\infty.$
\end{itemize}
\end{lemma}

\begin{lemma}\label{le:nehari}
For any $v\in E\setminus\{0\}$, there exists a unique $t(v)>0$ such that $t(v)v\in\mathcal{N}_\varepsilon$ and
\begin{equation*}
J_\varepsilon(t(v)v)=\displaystyle\max_{t\geq0}J_\varepsilon(tv).
\end{equation*}
\end{lemma}
{\bf Proof:} For any $v\in E\setminus\left\{0\right \}$, define $g(t)=J_\varepsilon(tv),\ t\in[0,+\infty).$ Then
\begin{equation*}
g(t)=\frac{t^2}{2}\|v\|^2_\varepsilon+\frac{bt^4}{4}(\int_{\mathbb{R}^3}|\nabla v|^2dx)^2-\frac{t^p}{p}\int_{\mathbb{R}^3}P(\varepsilon x)|v|^{p}dx-\frac{t^6}{6}\int_{\mathbb{R}^3}Q(\varepsilon x)|v|^{6}dx.
\end{equation*}
It is easy to see that $g(t)>0$ for $t>0$ small and $g(t)<0$ for $t>0$ large enough, so there exists $t_0>0$ such that
\begin{equation*}
g'(t_0)=0\quad \hbox{and}\quad g(t_0)=\max_{t\geq 0}g(t)=\max_{t\geq 0}J_\varepsilon(tv).
\end{equation*}
It follows from $g'(t_0)=0$ that $t_0v\in \mathcal{N}_\varepsilon$.
\par
If there exist $0<t_1<t_2$ such that $t_1v\in \mathcal{N}_\varepsilon$ and $t_2v\in \mathcal{N}_\varepsilon$. Then
\begin{equation*}
\frac{1}{t_1^2}\|v\|^2_\varepsilon+b(\int_{\mathbb{R}^3}|\nabla v|^2dx)^2=t_1^{p-4}\int_{\mathbb{R}^3}P(\varepsilon x)|v|^{p}dx+t_1^2\int_{\mathbb{R}^3}Q(\varepsilon x)|v|^{6}dx
\end{equation*}
and
\begin{equation*}
\frac{1}{t_2^2}\|v\|^2_\varepsilon+b(\int_{\mathbb{R}^3}|\nabla v|^2dx)^2=t_2^{p-4}\int_{\mathbb{R}^3}P(\varepsilon x)|v|^{p}dx+t_2^2\int_{\mathbb{R}^3}Q(\varepsilon x)|v|^{6}dx.
\end{equation*}
It follows that
\begin{equation*}
(\frac{1}{t_1^2}-\frac{1}{t_2^2})\|v\|^2_\varepsilon=(t_1^{p-4}-t_2^{p-4})\int_{\mathbb{R}^3}P(\varepsilon x)|v|^{p}dx+(t_1^2-t_2^2)\int_{\mathbb{R}^3}Q(\varepsilon x)|v|^{6}dx,
\end{equation*}
which is a contradiction.
\hfill{$\Box$}

\begin{lemma}\label{le:c-equal}
For any $\varepsilon>0$, let
\begin{equation*}
c_\varepsilon=\inf_{v\in \mathcal{N}_\varepsilon}J_\varepsilon(v),\quad c_\varepsilon^*=\inf_{v\in E\setminus\{0\}}\max_{t\geq0}J_\varepsilon(tv),
\quad c_\varepsilon^{**}=\inf_{\gamma\in \Gamma}\sup_{t\in [0,1]}J_\varepsilon(\gamma(t)),
\end{equation*}
where
\begin{equation*}
\Gamma=\{\gamma(t)\in C([0,1],E):\gamma(0)=0,J_\varepsilon(\gamma(1))<0\}.
\end{equation*}
Then, $c_\varepsilon=c_\varepsilon^*=c_\varepsilon^{**}.$
\end{lemma}
{\bf Proof:} We  divided the proof into three steps.
\par
{\bf Step1:} $c_\varepsilon^*=c_\varepsilon.$ By Lemma \ref{le:nehari}, we have
\begin{equation*}
c_\varepsilon^*=\displaystyle \inf_{v\in E\setminus\{0\}}\max_{t\geq0}J_\varepsilon(tv)
=\inf_{v\in E\setminus\{0\}}J_\varepsilon(t(v)v)=\inf_{v\in \mathcal{N}_\varepsilon}J_\varepsilon(v)=c_\varepsilon.
\end{equation*}
\par
{\bf Step2}. $c_\varepsilon^*\geq c_\varepsilon^{**}.$
From Lemma \ref{le:nehari}, for any $v\in E\setminus\{0\}$, there exists $T$ large enough, such that $J_\varepsilon(Tv)<0.$ Define $\gamma(t)=tTv$, $t\in[0,1]$, then
$\gamma(t)\in \Gamma$. Thus
$$c_\varepsilon^{**}=\inf_{\gamma\in \Gamma}\sup_{t\in [0,1]}J_\varepsilon(\gamma(t))\leq \sup_{t\in [0,1]}J_\varepsilon(\gamma(t))\leq \max_{t\geq0}J_\varepsilon(tv).$$
It follows that $c_\varepsilon^*\geq c_\varepsilon^{**}$.
\par
{\bf Step3}.  $ c_\varepsilon^{**}\geq c_\varepsilon.$
The manifold $\mathcal{N}_\varepsilon$ separates $E$ into two components. It is easy to know that the component containing the origin also contains a small ball around
the origin. It follows from  $g'(t)=\langle J'_\varepsilon(tv), v\rangle\geq 0$ for all $0\leq t\leq t(v)$ in Lemma \ref{le:nehari} that  $J_\varepsilon(v)\geq0$ in this component.
Thus every $\gamma\in \Gamma$ has to cross $\mathcal{N}_\varepsilon$. Then $ c_\varepsilon^{**}\geq c_\varepsilon$.
\hfill{$\Box$}

\begin{lemma}\label{le:psbounded} Any $(PS)_c$ sequence $\left\{v_n\right\}$
for $J_\varepsilon$ is bounded, and
\begin{equation*}
\limsup_{n\to \infty} \|v_n\|_\varepsilon\leq\sqrt{\frac{2p}{p-2}}c.
\end{equation*}
\end{lemma}
{\bf Proof:} Suppose that $\left\{v_n\right\}$ is a $(PS)_c$ sequence of
$J_\varepsilon,$ we have
\begin{equation*}
J_\varepsilon(v_n)\to c,\quad J'_\varepsilon(v_n)\to 0.
\end{equation*}
Thus
\begin{equation*}
\begin{aligned}
c+o(1)+o(1)\|v_n\|_\varepsilon
&=J_\varepsilon(v_n)-\displaystyle \frac{1}{p}\langle J'_\varepsilon(v_n),
v_n\rangle \\
&=\displaystyle(\frac{1}{2}-\frac{1}{p})\|v_n\|^2_\varepsilon+(\frac{1}{4}-\frac{1}{p})b(\int_{\mathbb{R}^3}|\nabla v_n|^2dx)^2\\
&\quad\displaystyle+(\frac{1}{p}-\frac{1}{6})\int_{\mathbb{R}^3}Q(\varepsilon x)v_n^6dx.
\end{aligned}
\end{equation*}
It follows that
\begin{equation*}
(\frac{1}{2}-\frac{1}{p})\|v_n\|^2_\varepsilon\leq c+o(1)+o(1)\|v_n\|_\varepsilon.
\end{equation*}
Then $\left\{v_n\right\}$ is bounded in $E,$ and
\begin{equation*}
\limsup_{n\to \infty} \|v_n\|_\varepsilon\leq\sqrt{\frac{2p}{p-2}}c.
\end{equation*}
\hfill{$\Box$}

\begin{lemma}\label{le:pssolution}
If $\{v_n\}$ is a $(PS)_{c_\varepsilon}$ sequence of $J_\varepsilon$ in $E$, then there exists $v\in E$ such that
$v_n\rightharpoonup v$ in  $E$ and $J'_\varepsilon(v)=0$.
\end{lemma}
{\bf Proof:} The proof is similar in \cite{Li-YeMMAS}, we give it for completeness. By Lemma \ref{le:psbounded}, we know that $\{v_n\}$ is bounded in $E$. Then,
up to a subsequence, we have
\begin{equation*}
\begin{aligned}
&v_n\rightharpoonup v\quad\hbox{in}\; E\\
&v_n\to v\quad\hbox{a.e. in}\  \mathbb{R}^3,\\
&v_n\rightharpoonup v
\quad\hbox{in}\  L^{q}(\mathbb{R}^3), \ \hbox{for }\ 2\leq q \leq 6.
\end{aligned}
\end{equation*}
For any function $\phi\in C_0^\infty(\mathbb{R}^3)$, since $J'_\varepsilon(v_n)\to 0$,  we have
\begin{equation*}
\begin{aligned}
o(1)&=\displaystyle\int_{\mathbb{R}^3}(a\nabla v_n \nabla \phi+ V(\varepsilon x)v_n \phi)dx+b\int_{\mathbb{R}^3}|\nabla v_n|^2dx\int_{\mathbb{R}^3}\nabla v_n \nabla \phi dx\\
&\quad\displaystyle-\int_{\mathbb{R}^3}P(\varepsilon x)|v_n|^{p-2}v_n\phi dx-\int_{\mathbb{R}^3}Q(\varepsilon x)|v_n|^{4}v_n\phi dx,
\end{aligned}
\end{equation*}
Assume that $\int_{\mathbb{R}^3}|\nabla v_n|^2dx\to A^2$, as $n\to \infty.$  Then, we have
\begin{equation*}
\begin{aligned}
0&=\displaystyle\int_{\mathbb{R}^3}(a\nabla v \nabla \phi+ V(\varepsilon x)v \phi)dx+bA^2\int_{\mathbb{R}^3}\nabla v \nabla \phi dx\\
&\quad\displaystyle-\int_{\mathbb{R}^3} P(\varepsilon x)|v|^{p-2}v\phi dx-\int_{\mathbb{R}^3}Q(\varepsilon x)|v|^{4}v\phi dx.
\end{aligned}
\end{equation*}
Thus, we can get
\begin{equation*}
\int_{\mathbb{R}^3}(a|\nabla v|^2+V(\varepsilon x)v^2)dx+ bA^2\int_{\mathbb{R}^3}|\nabla v|^2dx
=\int_{\mathbb{R}^3} P(\varepsilon x)|v|^{p}dx+\int_{\mathbb{R}^3}Q(\varepsilon x)|v|^{6}dx.
\end{equation*}
It is easy to know that $\int_{\mathbb{R}^3}|\nabla v|^2dx\leq A^2$. If $\int_{\mathbb{R}^3}|\nabla v|^2dx< A^2$, we have
\begin{equation*}
\int_{\mathbb{R}^3}(a|\nabla v|^2+V(\varepsilon x)v^2)dx+ b(\int_{\mathbb{R}^3}|\nabla v|^2dx)^2
<\int_{\mathbb{R}^3} P(\varepsilon x)|v|^{p}dx+\int_{\mathbb{R}^3}Q(\varepsilon x)|v|^{6}dx.
\end{equation*}
Then, there exists $0<t_0<1$ such that $t_0v\in \mathcal{N}_\varepsilon$.  Thus,  we have
\begin{equation*}
\begin{aligned}
c_\varepsilon&\leq J_\varepsilon(t_0v)-\displaystyle\frac{1}{4}\langle J_\varepsilon'(t_0v),t_0v\rangle\\
&=\displaystyle\frac{t_0^2}{4}\int_{\mathbb{R}^3}(a|\nabla v|^2+V(\varepsilon x)v^2)dx+(\frac{t_0^p}{4}-\frac{t_0^p}{p})\int_{\mathbb{R}^3} P(\varepsilon x)|v|^{p}dx+(\frac{t_0^6}{4}-\frac{t_0^6}{6})\int_{\mathbb{R}^3}Q(\varepsilon x)|v|^{6}dx\\
&<\displaystyle\frac{1}{4}\int_{\mathbb{R}^3}(a|\nabla v|^2+V(\varepsilon x)v^2)dx+(\frac{1}{4}-\frac{1}{p})\int_{\mathbb{R}^3} P(\varepsilon x)|v|^{p}dx+(\frac{1}{4}-\frac{1}{6})\int_{\mathbb{R}^3}Q(\varepsilon x)|v|^{6}dx\\
&\leq\displaystyle\liminf_{n\to\infty}\Big[\frac{1}{4}\int_{\mathbb{R}^3}(a|\nabla v_n|^2+V(\varepsilon x)v_n^2)dx+(\frac{1}{4}-\frac{1}{p})\widetilde\int_{\mathbb{R}^3} P(\varepsilon x)|v_n|^{p}dx\\
&\quad\displaystyle+(\frac{1}{4}-\frac{1}{6})\int_{\mathbb{R}^3}Q(\varepsilon x)|v_n|^{6}dx\Big]\\
&=\displaystyle\liminf_{n\to\infty}\Big[J_{\varepsilon}(v_n)-\displaystyle\frac{1}{4}\langle J_{\varepsilon}'(v_n),v_n\rangle\Big]\\
&=c_{\varepsilon}.
\end{aligned}
\end{equation*}
Then, as $n\to \infty$, we have $\int_{\mathbb{R}^3}|\nabla v_n|^2dx\to \int_{\mathbb{R}^3}|\nabla v|^2dx$. Therefore, $J'_\varepsilon(v)=0$.
\hfill{$\Box$}
\vspace{6mm}

In order to investigate \eqref{transeq1}, we need some results about \eqref{transeq1} with constant coefficients. Consider the following problem
\begin{equation}\label{transeqc1}
-(a+ b\int_{\mathbb{R}^3}|\nabla v|^2dx)\Delta v+kv=\tau|v|^{p-2}v+\nu|v|^4v,\quad x\in\mathbb{R}^3.
\end{equation}
where $k$, $\tau$ and $\nu$ are positive constants. The associated energy functional is
\begin{equation*}
\Phi^\ast_{k\tau\nu}(v)=\frac{1}{2}\|v\|^2_k+\frac{b}{4}(\int_{\mathbb{R}^3}|\nabla v|^2dx)^2-\frac{\tau}{p}\int_{\mathbb{R}^3}|v|^{p}dx-\frac{\nu}{6}\int_{\mathbb{R}^3}|v|^{6}dx,
\end{equation*}
where $\|v\|^2_k:=\int_{\mathbb{R}^3}(a|\nabla v|^2+kv^2)dx$ and
$\|v\|_k$ is an equivalent norm on $E$.

By Lemma \ref{le:c-equal}, we have $$m^*_{k\tau\nu}:=\inf_{v\in E\setminus\{0\}}\max_{t\geq0}\Phi^*_{k\tau\nu}(tv)=\inf_{v\in \mathcal{N}^\ast_{k\tau\nu}}\Phi^*_{k\tau\nu}(v),$$ where $\mathcal{N}^\ast_{k\tau\nu}=\left\{v\in E\setminus \left\{0\right\}| \langle (\Phi^*_{k\tau\nu})'(v), v\rangle=0\right\}.$

\begin{lemma}\cite{MR3347410}\label{inequality}
For $t,\ s>0$ and $\lambda$ is a positive constant, the following system
\begin{equation*}
\begin{cases}
\Phi(t,s)=t-aS\lambda^{\frac{-1}{3}}(t+s)^{\frac{1}{3}}=0, \\
\Psi(t,s)=s-bS^2\lambda^{\frac{-2}{3}}(t+s)^{\frac{2}{3}}=0.\\
\end{cases}
\end{equation*}
has a unique solution $(t_0,s_0)$. Moreover, if $\Phi(t,s)\geq0$ and $\Psi(t,s)\geq0$, then $t\geq t_0$, $s\geq s_0$.
\end{lemma}
\par
In order to establish the existence result, we need to show that the mountain pass value is less than the critical level. This result can be found in the following Lemma.
\begin{lemma}\label{le:c-level}\cite{MR3347410}
For any $\varepsilon>0$ and $Q(x)\equiv q>0,$ we have $c_\varepsilon<\displaystyle\frac{ab}{4q}S^3+\frac{(b^2S^4+4qaS)^{3/2}}{24q^2}+\frac{b^3S^6}{24q^2}$.
\end{lemma}

\begin{lemma}\label{le:constant-achieved}
Problem \eqref{transeqc1} has at least one positive ground state solution.
\end{lemma}
{\bf Proof:} By Lemma \ref{le:moutain-pass geometry} and Lemma \ref{le:c-equal}, there exits a sequence $\left\{v_n\right\}$ which is a $(PS)_{m^*_{k\tau\nu}}$ sequence of
$\Phi^\ast_{k\tau\nu}(v)$. From Lemma \ref{le:psbounded}, we know that $\{v_n\}$ is bounded in $E$. Hence, up to a
subsequence, we have
\begin{equation*}
\begin{aligned}
&v_n\rightharpoonup v\quad\hbox{in}\; E\\
&v_n\to v\quad\hbox{a.e. in}\  \mathbb{R}^3,\\
&v_n\rightharpoonup v
\quad\hbox{in}\  L^{q}(\mathbb{R}^3), \ \hbox{for }\ 2\leq q \leq 6.
\end{aligned}
\end{equation*}
It follows from Lemma \ref{le:pssolution} that $(\Phi^{\ast}_{k\tau\nu})'(v)=0$ .
\par
Since $\{v_n\}$ is a $(PS)_{m^*_{k\tau\nu}}$ sequence of $\Phi^{\ast}_{k\tau\nu}$, we have
\begin{equation*}
\begin{aligned}
o(1)&=\displaystyle\langle (\Phi^{\ast}_{k\tau\nu})'(v_n),v_n\rangle\\
&=\displaystyle\|v_n\|^2_k+b(\int_{\mathbb{R}^3}|\nabla v_n|^2dx)^2-\tau\int_{\mathbb{R}^3}|v_n|^{p}dx-\nu\int_{\mathbb{R}^3}|v_n|^{6}dx.
\end{aligned}
\end{equation*}
Since $\{v_n\}$ is bounded in $E$, as $n\rightarrow\infty$, we can assume
\begin{equation}\label{lim1}
\|v_n\|^2_k\rightarrow l_1,
\end{equation}
\begin{equation}\label{lim2}
b(\int_{\mathbb{R}^3}|\nabla v_n|^2dx)^2\rightarrow l_2,
\end{equation}
and
\begin{equation}\label{lim3}
\tau\int_{\mathbb{R}^3}|v_n|^{p}dx+\nu\int_{\mathbb{R}^3}|v_n|^{6}dx\rightarrow l_3.
\end{equation}
Then we have $l_3=l_1+l_2$.

If $l_3=0$, we have $v_n\rightarrow 0$ in $E.$ Then $\Phi^{\ast}_{k\tau\nu}(v_n)\rightarrow 0,$ which contradicts $m^*_{k\tau\nu}>0.$ Thus, $l_3\neq0$.

In \eqref{lim3}, if $\displaystyle\int_{\mathbb{R}^3}|v_n|^{p}dx\rightarrow 0,$ then we have
\begin{equation}\label{lim4}
\nu\int_{\mathbb{R}^3}|v_n|^6dx\rightarrow l_3.
\end{equation}
By the definition of the best constant $S$, we have
\begin{equation*}
a^3\nu\int_{\mathbb{R}^3}|v_n|^6dx\leq a^3\nu(S^{-1}\int_{\mathbb{R}^3}|\nabla v_n|^2dx)^3
\leq\nu S^{-3}\|v_n\|^6_k
\end{equation*}
and
\begin{equation*}
b(\nu\int_{\mathbb{R}^3}|v_n|^6dx)^{2/3}\leq b\nu^{2/3}(S^{-1}\int_{\mathbb{R}^3}|\nabla v_n|^2dx)^2\leq  b\nu^{2/3}S^{-2}(\int_{\mathbb{R}^3}|\nabla v_n|^2dx)^2.
\end{equation*}
Letting $n\to\infty$ in the above two inequalities, we have
\begin{equation*}
a^3l_3\leq \nu S^{-3}l^3_1\quad  \hbox{and} \quad  bl_3^{2/3}\leq \nu^{2/3}S^{-2}l_2.
\end{equation*}
Therefore,
\begin{equation*}
l_1\geq aS\nu^{-\frac{1}{3}}(l_1+l_2)^{\frac{1}{3}}\quad  \hbox{and} \quad  l_2\geq bS^2\nu^{-\frac{2}{3}}(l_1+l_2)^{\frac{2}{3}}.
\end{equation*}
By Lemma \ref{inequality}, we have
\begin{equation*}
\frac{1}{3}l_1+\frac{1}{12}l_2\geq c^*:=\frac{ab}{4\nu}S^3+\frac{(b^2S^4+4\nu aS)^{3/2}}{24\nu^2}+\frac{b^3S^6}{24\nu^2}.
\end{equation*}
On the other hand,
\begin{equation*}
\begin{aligned}
m^*_{k\tau\nu}&=\Phi^\ast_{k\tau\nu}(v_n)+o(1)\\
&=\displaystyle\frac{1}{2}\|v_n\|^2_k+\frac{b}{4}(\int_{\mathbb{R}^3}|\nabla v_n|^2dx)^2-\frac{\tau}{p}\int_{\mathbb{R}^3}|v_n|^{p}dx-\frac{\nu}{6}\int_{\mathbb{R}^3}|v_n|^{6}dx+o(1)\\
&=\displaystyle\frac{1}{2}\|v_n\|^2_k+\frac{b}{4}(\int_{\mathbb{R}^3}|\nabla v_n|^2dx)^2-\frac{\nu}{6}\int_{\mathbb{R}^3}|v_n|^{6}dx+o(1)\\
&=\displaystyle\frac{1}{2}\|v_n\|^2_k+\frac{b}{4}(\int_{\mathbb{R}^3}|\nabla v_n|^2dx)^2-\frac{1}{6}(\|v_n\|^2_k+b(\int_{\mathbb{R}^3}|\nabla v_n|^2dx)^2)+o(1)\\
&=\displaystyle\frac{1}{3}\|v_n\|^2_k+\frac{b}{12}(\int_{\mathbb{R}^3}|\nabla v_n|^2dx)^2+o(1)\\
\end{aligned}
\end{equation*}
Using \eqref{lim1} and \eqref{lim2} in the above expression, we obtain  $$m^*_{k\tau\nu}= \frac{1}{3}l_1+\frac{1}{12}l_2\geq c^*,$$
which contradicts Lemma \ref{le:c-level}.
\par
Therefore, $\displaystyle\int_{\mathbb{R}^3}|v_n|^{p}\rightarrow l_4>0,$ as $n\rightarrow\infty.$ Then, by Lions's Lemma,
there exists $(y_n)\subset \mathbb{R}^3$, $\rho,\eta>0$ such that
\begin{equation}\label{lionslemma}
\displaystyle\limsup_{n\to\infty}\int_{B_\rho(y_n)}|v_n|^2dx\geq \eta.
\end{equation}
Let $\widetilde{v}_n(x)=v_n(x+y_n)$, then $\|\widetilde{v}_n\|_k\leq C$ in $E$, so there exists $\widetilde{v}\in E$ such that
$\widetilde{v}_n\rightharpoonup \widetilde{v}$ in $E$ and $\widetilde{v}_n\rightarrow \widetilde{v}$ a.e in $\mathbb{R}^3$, by \eqref{lionslemma}, we get $\widetilde{v}\neq0$.
\par
It is easy to prove that
\begin{equation*}
(\Phi^{\ast}_{k\tau\nu})(\widetilde{v}_n)\to m^*_{k\tau\nu},\quad (\Phi^{\ast}_{k\tau\nu})'(\widetilde{v}_n)\to 0.
\end{equation*}
Then, we have $(\Phi^{\ast}_{k\tau\nu})'(\widetilde{v})=0$ and $\widetilde{v}\in \mathcal{N}^\ast_{k\tau\nu}$.
\par
Moreover,
\begin{equation*}
\begin{aligned}
m^*_{k\tau\nu}&=\displaystyle\lim_{n\to\infty}\Big[(\Phi^{\ast}_{k\tau\nu})(\widetilde{v}_n)-\displaystyle\frac{1}{4}\langle (\Phi^{\ast}_{k\tau\nu})'(\widetilde{v}_n),\widetilde{v}_n\rangle\Big]\\
&=\displaystyle\lim_{n\to\infty}\Big[\frac{1}{4}\|\widetilde{v}_n\|^2_k+
(\frac{1}{4}-\frac{1}{p})\tau\int_{\mathbb{R}^3}|\widetilde{v}_n|^{p}dx+(\frac{1}{4}-\frac{1}{6})\nu\int_{\mathbb{R}^3}|\widetilde{v}_n|^{6}dx\Big]\\
&\geq\displaystyle\frac{1}{4}\|\widetilde{v}\|^2_k+
(\frac{1}{4}-\frac{1}{p})\tau\int_{\mathbb{R}^3}|\widetilde{v}|^{p}dx+(\frac{1}{4}-\frac{1}{6})\nu\int_{\mathbb{R}^3}|\widetilde{v}|^{6}dx\\\
&=(\Phi^{\ast}_{k\tau\nu})(\widetilde{v})-\displaystyle\frac{1}{4}\langle (\Phi^{\ast}_{k\tau\nu})'(\widetilde{v}),\widetilde{v}\rangle\\
&=(\Phi^{\ast}_{k\tau\nu})(\widetilde{v}),
\end{aligned}
\end{equation*}
which means  $(\Phi^{\ast}_{k\tau\nu})(\widetilde{v})=m^*_{k\tau\nu}.$ It is easy to know that $|\widetilde{v}|\in N^\ast_{k\tau\nu}$ and $(\Phi^{\ast}_{k\tau\nu})(|\widetilde{v}|)=m^*_{k\tau\nu}.$ According to the proof of
Theorem 4.3 in \cite{Willem1996}, we can show that  $(\Phi^{\ast}_{k\tau\nu})'(|\widetilde{v}|)=0.$ Without loss of generality, we can assume $\widetilde{v}\geq0.$ By the theory of
elliptic regularity, $\widetilde{v}\in C^2(\mathbb{R}^3)$, and by using strong maximum principle, we get $\widetilde{v}>0$ in $\mathbb{R}^3$.
\hfill{$\Box$}

\begin{lemma}\label{le:comparel} For the positive constants $k_i$, $\tau_i$ and $\nu_i$, $i=1,2$.  If $$\min\left\{k_2-k_1,\tau_1-\tau_2,\nu_1-\nu_2\right\}\geq 0,$$ then $m^*_{k_1\tau_1\nu_1}\leq m^*_{k_2\tau_2\nu_2}$. Additionally, if $\max\left\{k_2-k_1,\tau_1-\tau_2,\nu_1-\nu_2\right\}> 0,$ then $m^*_{k_1\tau_1\nu_1}< m^*_{k_2\tau_2\nu_2}$.
\end{lemma}
{\bf Proof:} By Lemma \ref{le:constant-achieved}, there exist $v\in E$ such that $\Phi^*_{k_2\tau_2\nu_2}(v)= m^*_{k_2\tau_2\nu_2}=\displaystyle\max_{t\geq0}\Phi^*_{k_2\tau_2\nu_2}(tv).$
By Lemma \ref{le:nehari}, there exists $t_0>0$ such that $\Phi^*_{k_1\tau_1\nu_1}(t_0v)=\displaystyle\max_{t\geq0}\Phi^*_{k_1\tau_1\nu_1}(tv).$ Then
$$m^*_{k_1\tau_1\nu_1}\leq\displaystyle\max_{t\geq0}\Phi^*_{k_1\tau_1\nu_1}(tv)=\Phi^*_{k_1\tau_1\nu_1}(t_0v)\leq\Phi^*_{k_2\tau_2\nu_2}(t_0v)\leq\Phi^*_{k_2\tau_2\nu_2}(v)=m^*_{k_2\tau_2\nu_2}.$$
\hfill{$\Box$}

\section{An Auxiliary Problem}
In this section, we mainly consider an auxiliary problem, for $c\in[V_{min},V_{max}],$ $d\in[P_{min},P_{max}]$ and $e \in[Q_{min},Q_{max}],$ we define
\begin{equation*}
\begin{aligned}
V^c_\varepsilon(x)=\max\{c,V(\varepsilon x)\},\\
P^d_\varepsilon(x)=\min\{d,P(\varepsilon x)\},\\
Q^e_\varepsilon(x)=\min\{e,Q(\varepsilon x)\}.
\end{aligned}
\end{equation*}
When $x=0$, we set $V^c(0)=\max\{c,V(0)\},$ $P^d(0)=\min\{d,P(0)\}$ and
$Q^e(0)=\min\{e,Q(0)\}.$
\par
Consider the following equation
\begin{equation}\label{tmaineq1}
-(a+b\int_{\mathbb{R}^3}|\nabla v|^2dx)\Delta v+V^c_\varepsilon(x)v=P^d_\varepsilon(x)|v|^{p-2}v+Q^e_\varepsilon(x)|v|^4v,\quad x\in\mathbb{R}^3.
\end{equation}
whose energy functional is
\begin{equation*}
\begin{aligned}
J^{cde}_\varepsilon(v)=&\displaystyle\frac{1}{2}\int_{\mathbb{R}^3}(|\nabla
v|^2+V^c_\varepsilon(x)v^2)dx+\frac{b}{4}(\int_{\mathbb{R}^3}|\nabla v|^2dx)^2\\
&\quad-\displaystyle\frac{1}{p}\int_{\mathbb{R}^3}P^d_\varepsilon(x)|v|^{p}dx-\frac{1}{6}\int_{\mathbb{R}^3}Q^e_\varepsilon(x)|v|^6dx.
\end{aligned}
\end{equation*}
By Lemma \ref{le:c-equal}, we have
\begin{equation*}
c^{cde}_\varepsilon:=\inf_{v\in E\setminus\{0\}}\max_{t\geq0}J^{cde}_\varepsilon(tv)=\inf_{v\in \mathcal{N}^{cde}_\varepsilon}J^{cde}_\varepsilon(v),
\end{equation*}
where  $\mathcal{N}^{cde}_\varepsilon=\left\{v\in E\setminus \left\{0\right\}| \langle (J^{cde}_\varepsilon)'(v), v\rangle=0\right\}.$

\begin{lemma}\label{cv(y)}
For any $y\in \mathbb{R}^3$, $\displaystyle\limsup_{\varepsilon\rightarrow 0}c_\varepsilon\leq m^*_{V(y)P(y)Q(y)}.$
\end{lemma}
{\bf Proof:} Let $v$ be a ground state solution of
\begin{equation*}
-(a+b\int_{\mathbb{R}^3}|\nabla v|^2dx)\Delta v+V(y)v=P(y)|v|^{p-2}v+Q(y)|v|^4v,\quad x\in\mathbb{R}^3.
\end{equation*}
Then, we have $v\in \mathcal{N}^\ast_{V(y)P(y)Q(y)}$ and $\Phi^*_{V(y)P(y)Q(y)}(v)= m^*_{V(y)P(y)Q(y)}.$
\par
Define $v_\varepsilon(x)=v(x-\frac{y}{\varepsilon})$. From Lemma \ref{le:nehari}, we know that there exists a unique $t_\varepsilon>0$ satisfying
\begin{equation*}
t_\varepsilon v_\varepsilon \in \mathcal{N}_\varepsilon \quad \hbox{and} \quad J_\varepsilon(t_\varepsilon v_\varepsilon)=\displaystyle\sup_{t\geq0}J_\varepsilon(tv_\varepsilon).
\end{equation*}
By $t_\varepsilon v_\varepsilon \in \mathcal{N}_\varepsilon$, we can prove $t_\varepsilon$ is bounded. It follows from Lemma \ref{le:lowboundedinN} that $t_\varepsilon$ has a positive lower bound. Without loss of generality, we assume $t_\varepsilon\rightarrow t_0>0.$
\par
Since $t_\varepsilon v_\varepsilon \in \mathcal{N}_\varepsilon$, we have
\begin{equation*}
\begin{aligned}
&\displaystyle t_\varepsilon^2\int_{\mathbb{R}^3}(a|\nabla
v_\varepsilon|^2+V(\varepsilon x)v_\varepsilon^2)dx+bt_\varepsilon^4(\int_{\mathbb{R}^3}|\nabla v_\varepsilon|^2dx)^2\\
&\displaystyle=t_\varepsilon^p\int_{\mathbb{R}^3}P(\varepsilon x)|v_\varepsilon|^{p}dx+t_\varepsilon^6\int_{\mathbb{R}^3}Q(\varepsilon x)|v_\varepsilon|^6dx.
\end{aligned}
\end{equation*}
Thus, we have
\begin{equation*}
\begin{aligned}
&\displaystyle t_\varepsilon^2\int_{\mathbb{R}^3}(a|\nabla
v|^2+V(\varepsilon x+y)v^2)dx+bt_\varepsilon^4(\int_{\mathbb{R}^3}|\nabla v|^2dx)^2\\
&\displaystyle=t_\varepsilon^p\int_{\mathbb{R}^3}P(\varepsilon x+y)|v|^{p}dx+t_\varepsilon^6\int_{\mathbb{R}^3}Q(\varepsilon x+y)|v|^6dx.
\end{aligned}
\end{equation*}
Letting $\varepsilon\rightarrow 0,$ by Lebesgue Dominated Convergence Theorem, we get
\begin{equation*}
\begin{aligned}
&\displaystyle t_0^2\int_{\mathbb{R}^3}(a|\nabla
v|^2+V(y)v^2)dx+bt_0^4(\int_{\mathbb{R}^3}|\nabla v|^2dx)^2\\
&\displaystyle=t_0^p\int_{\mathbb{R}^3}P(y)|v|^{p}dx+t_0^6\int_{\mathbb{R}^3}Q(y)|v|^6dx,
\end{aligned}
\end{equation*}
which means $t_0v\in \mathcal{N}^\ast_{V(y)P(y)Q(y)}.$ By using Lemma \ref{le:nehari} and $v\in \mathcal{N}^\ast_{V(y)P(y)Q(y)}$, we obtain $t_0=1.$
So we get
\begin{equation*}
\begin{aligned}
c_\varepsilon
&\leq J_\varepsilon(t_\varepsilon v_\varepsilon)\\
&=\displaystyle\frac{t_\varepsilon^2}{2}\int_{\mathbb{R}^3}( a|\nabla v_\varepsilon|^2+V(\varepsilon x)  v_\varepsilon^2)dx+\frac{bt_\varepsilon^4}{4}(\int_{\mathbb{R}^3}|\nabla v_\varepsilon|^2dx)^2
-\displaystyle\frac{t_\varepsilon^p}{p}\int_{\mathbb{R}^3}P(\varepsilon x)|v_\varepsilon|^{p}dx\\
&\quad-\displaystyle\frac{t_\varepsilon^6}{6}\int_{\mathbb{R}^3}Q(\varepsilon x)|v_\varepsilon|^6dx\\
&=\displaystyle\frac{t_\varepsilon^2}{2}\int_{\mathbb{R}^3}( a|\nabla v|^2+V(\varepsilon x+y)  v^2)dx+\frac{bt_\varepsilon^4}{4}(\int_{\mathbb{R}^3}|\nabla v|^2dx)^2
-\displaystyle\frac{t_\varepsilon^p}{p}\int_{\mathbb{R}^3}P(\varepsilon x+y)|v|^{p}dx\\
&\quad-\displaystyle\frac{t_\varepsilon^6}{6}\int_{\mathbb{R}^3}Q(\varepsilon x+y)|v|^6dx\\
&=\displaystyle\frac{1}{2}\int_{\mathbb{R}^3}( a|\nabla v|^2+V(y)v^2)dx+\frac{b}{4}(\int_{\mathbb{R}^3}|\nabla v|^2dx)^2
-\displaystyle\frac{1}{p}\int_{\mathbb{R}^3}P(y)|v|^{p}dx\\
&\quad-\displaystyle\frac{1}{6}\int_{\mathbb{R}^3}Q(y)|v|^6dx+o_\varepsilon(1)\\
&=\Phi^*_{V(y)P(y)Q(y)}(v)+ o_\varepsilon(1)\\
&=m^*_{V(y)P(y)Q(y)}+o_\varepsilon(1).
\end{aligned}
\end{equation*}
Therefore, $\displaystyle\limsup_{\varepsilon\rightarrow 0}c_\varepsilon\leq m^*_{V(y)P(y)Q(y)}.$
\hfill{$\Box$}

\begin{lemma}\label{cbigm} For any $\varepsilon>0$, we have $c^{V(0) dQ_{max}}_\varepsilon\geq m^*_{V(0) dQ_{max}}.$\end{lemma}
{\bf Proof:} For any $v\in E,$ we have $J^{V(0)dQ_{max}}_\varepsilon(tv)\geq \Phi^\ast_{V(0) dQ_{max}}(tv).$ So
\begin{equation*}
\inf_{u\in E}\max_{t>0}J^{V(0)dQ_{max}}_\varepsilon(tv)\geq \inf_{u\in E}\max_{t>0}\Phi^\ast_{V(0)dQ_{max}}(tv).
\end{equation*}
By Lemma \ref{le:c-equal}, the proof is completed.
\hfill{$\Box$}

\section{Proof of the Main Results}

\begin{lemma} \label{lemma4.1}
Suppose that the potential functions $V(x)$, $P(x)$ and $Q(x)$ satisfy conditions $(PQ1)$ and $(PQ2)$. Then for any $\varepsilon>0$ small enough, problem \eqref{maineq1} has at least one positive ground state solution.
\end{lemma}
{\bf Proof:} By Lemma \ref{le:c-equal}, we can choose a sequence $\left\{v_n\right\}\subset \mathcal{N}_\varepsilon$ such that $J_\varepsilon(v_n)\rightarrow c_\varepsilon.$ In view of Ekeland's variational principle, the sequence can be chosen to be a $(PS)_{c_\varepsilon}$ sequence of
$J_\varepsilon(v)$. From Lemma \ref{le:psbounded}, we know that $\left\{v_n\right\}$ is bounded in $E$. Hence, up to a
subsequence, we have
\begin{equation*}
\begin{aligned}
&v_n\rightharpoonup v\quad\hbox{in}\; E\\
&v_n\to v\quad\hbox{a.e. in}\  \mathbb{R}^3,\\
&v_n\rightharpoonup v
\quad\hbox{in}\  L^{q}(\mathbb{R}^3), \ \hbox{for }\ 2\leq q \leq 6.
\end{aligned}
\end{equation*}
It follows from Lemma \ref{le:pssolution} that $J_{\varepsilon}'(v_\varepsilon)=0$. Now we prove $v_\varepsilon\neq0.$
\par
\textbf{Claim 1:} there exist $(y_n)\subset \mathbb{R}^3$ and $\rho,\eta>0$ such that
\begin{equation}\label{lionslemma3}
\displaystyle\limsup_{n\to\infty}\int_{B_\rho(y_n)}|v_n|^2dx\geq \eta.
\end{equation}
Otherwise, we have $v_n\rightarrow0$ in $L^q(\mathbb{R}^3), \ 2<q<6.$
\par
Since $\left\{v_n\right\}$ is a $(PS)_{c_\varepsilon}$ sequence of $J_\varepsilon$, we have
\begin{equation*}
\begin{aligned}
o(1)&=\displaystyle\langle J_\varepsilon'(v_n),v_n\rangle\\
&=\displaystyle\|v_n\|^2_\varepsilon+b(\int_{\mathbb{R}^3}|\nabla v_n|^2dx)^2-\int_{\mathbb{R}^3}P(\varepsilon x)|v_n|^{p}dx-\int_{\mathbb{R}^3}Q(\varepsilon x)|v_n|^{6}dx.
\end{aligned}
\end{equation*}
Thus,
\begin{equation*}
o(1)=\|v_n\|^2_\varepsilon+b(\int_{\mathbb{R}^3}|\nabla v_n|^2dx)^2-\int_{\mathbb{R}^3}Q(\varepsilon x)|v_n|^{6}dx.
\end{equation*}
Since $\{v_n\}$ is bounded in $E$, as $n\rightarrow\infty$, we can assume
\begin{equation}\label{lim12}
\|v_n\|^2_\varepsilon\rightarrow l_1,
\end{equation}
\begin{equation}\label{lim22}
b(\int_{\mathbb{R}^3}|\nabla v_n|^2dx)^2\rightarrow l_2,
\end{equation}
and
\begin{equation}\label{lim32}
\int_{\mathbb{R}^3}Q(\varepsilon x)|v_n|^{6}dx\rightarrow l_3.
\end{equation}
Then we have $l_3=l_1+l_2$.
\par
If $l_3=0$, we have $v_n\rightarrow 0$ in $E.$ Then $J_\varepsilon(v_n)\rightarrow 0,$ which contradicts $c_\varepsilon>0.$ Thus, $l_3\neq0$.
By the definition of the best constant $S$, we have
\begin{equation*}
S\leq \frac{\displaystyle\frac{1}{a}\int_{\mathbb{R}^3}a|\nabla v_n|^2dx}{\displaystyle(Q^{-1}_{max}\int_{\mathbb{R}^3}Q(\varepsilon x)|v_n|^{6}dx)^{1/3}}
\leq \frac{\displaystyle\frac{1}{a}\|v_n\|^2_\varepsilon}{\displaystyle(Q^{-1}_{max}\int_{\mathbb{R}^3}Q(\varepsilon x)|v_n|^{6}dx)^{1/3}}
\end{equation*}
and
\begin{equation*}
bS^2\leq \frac{ \displaystyle b(\int_{\mathbb{R}^3}|\nabla v_n|^2dx)^2}{\displaystyle Q^{-2/3}_{max}(\int_{\mathbb{R}^3}Q(\varepsilon x)|v_n|^{6}dx)^{2/3}}.
\end{equation*}
Letting $n\to\infty$ in the above two inequalities, we have
\begin{equation*}
aSl_3^{1/3}\leq Q_{max}^{1/3}l_1 \quad  \hbox{and} \quad  bS^2l_3^{2/3}\leq Q_{max}^{2/3}l_2.
\end{equation*}
Therefore,
\begin{equation*}
l_1\geq aSQ_{max}^{-\frac{1}{3}}(l_1+l_2)^{\frac{1}{3}}\quad  \hbox{and} \quad  l_2\geq bS^2Q_{max}^{-\frac{2}{3}}(l_1+l_2)^{\frac{2}{3}}.
\end{equation*}
By Lemma \ref{inequality}, we have
\begin{equation*}
\frac{1}{3}l_1+\frac{1}{12}l_2\geq c^*:=\frac{ab}{4Q_{max}}S^3+\frac{(b^2S^4+4Q_{max} aS)^{3/2}}{24Q_{max}^2}+\frac{b^3S^6}{24Q_{max}^2}.
\end{equation*}
On the other hand,
\begin{equation*}
\begin{aligned}
c_\varepsilon&=J_\varepsilon(v_n)+o(1)\\
&=\displaystyle\frac{1}{2}\|v_n\|^2_\varepsilon+\frac{b}{4}(\int_{\mathbb{R}^3}|\nabla v_n|^2dx)^2-\frac{1}{p}\int_{\mathbb{R}^3}P(\varepsilon x)|v_n|^{p}dx-\frac{1}{6}\int_{\mathbb{R}^3}Q(\varepsilon x)|v_n|^{6}dx+o(1)\\
&=\displaystyle\frac{1}{2}\|v_n\|^2_\varepsilon+\frac{b}{4}(\int_{\mathbb{R}^3}|\nabla v_n|^2dx)^2-\frac{1}{6}\int_{\mathbb{R}^3}Q(\varepsilon x)|v_n|^{6}dx+o(1)\\
&=\displaystyle\frac{1}{2}\|v_n\|^2_\varepsilon+\frac{b}{4}(\int_{\mathbb{R}^3}|\nabla v_n|^2dx)^2-\frac{1}{6}(\|v_n\|^2_\varepsilon+b(\int_{\mathbb{R}^3}|\nabla v_n|^2dx)^2)+o(1)\\
&=\displaystyle\frac{1}{3}\|v_n\|^2_\varepsilon+\frac{b}{12}(\int_{\mathbb{R}^3}|\nabla v_n|^2dx)^2+o(1).\\
\end{aligned}
\end{equation*}
Using \eqref{lim12} and \eqref{lim22} in the above expression, we obtain
\begin{equation*}
c_\varepsilon= \frac{1}{3}l_1+\frac{1}{12}l_2\geq c^*,
\end{equation*}
which contradicts Lemma \ref{le:c-level}. Therefore, \textbf{Claim 1} holds.
\par
Let $\widetilde{v}_n(x)=v_n(x+y_n)$. Then $\|\widetilde{v}_n\|_\varepsilon\leq C$ in $E$. So there exists $\widetilde{v}\in E$ such that
\begin{equation*}
\begin{aligned}
&\widetilde{v}_n\rightharpoonup \widetilde{v}\quad\hbox{in}\; E\\
&\widetilde{v}_n\rightarrow \widetilde{v}\quad\hbox{a.e. in}\  \mathbb{R}^3.\\
\end{aligned}
\end{equation*}
By \eqref{lionslemma3}, we get $\widetilde{v}\neq0$. Thus, there exists $\delta>0$ satisfying $\mu\{x:|\widetilde{v}(x)|>\delta\}>0.$ By $(PQ2)$, without loss of generality, we may assume $x^*=0\in \mathcal{P\cap Q},$ such that $\beta:=V(0)\leq V(x)$ for $|x|\geq R.$ Let $P_{\infty}< d<P_{max}.$ For $v_n$, there exists $t_n>0$ such that $t_nv_n\in N_\varepsilon^{\beta dQ_{max}}$.
\par
\textbf{Claim 2:} $t_n$ is bounded.
\par
From $t_nv_n\in \mathcal{N}_\varepsilon^{\beta dQ_{max}}$, we get
\begin{equation*}
\begin{aligned}
&\displaystyle \int_{\mathbb{R}^3}( a|\nabla v_n|^2+V^\beta_\varepsilon(x) v_n^2) dx+bt_n^2(\int_{\mathbb{R}^3}|\nabla v_n|^2dx)^2\\
=&\displaystyle t_n^{p-2}\int_{\mathbb{R}^3}P^d_\varepsilon(x)|v_n|^pdx+ t_n^4\int_{\mathbb{R}^3}Q(\varepsilon x)| v_n|^6dx.
\end{aligned}
\end{equation*}
Then, by the boundness of $\{v_n\}$ in $E$, we have
\begin{equation*}
C+t_n^2C\geq t_n^4\int_{\mathbb{R}^3}Q(\varepsilon x)| v_n|^6dx\geq C\displaystyle t_n^4\int_{\{x:|\widetilde{v}(x)|>\delta\}}|\widetilde{v}_n|^6dx.
\end{equation*}
By Egoroff Theorem, there exists $E_0\subset\{x:|\widetilde{v}(x)|>\delta\}$  such that  $\mu\{\{x:|\widetilde{v}(x)|>\delta\}\setminus E_0\}>0$ and $\widetilde{v}_n\rightarrow \widetilde{v}$ uniformly in $\{x:|\widetilde{v}(x)|>\delta\}\setminus E_0$. Thus, we can obtain
\begin{equation*}
C+t_n^2C\geq C\displaystyle t_n^4\int_{\{x:|\widetilde{v}(x)|>\delta\}\setminus E_0}|\widetilde{v}_n|^6dx\geq Ct_n^4.
\end{equation*}
It follows that \textbf{Claim 2} holds.
\par
\textbf{Claim 3:}$J^{\beta dQ_{max}}_\varepsilon(t_nv_n)=J_\varepsilon(t_nv_n)+o_n(1).$
\par
First, we note
\begin{equation*}
\begin{aligned}
J^{\beta dQ_{max}}_\varepsilon(t_nv_n)=&J_\varepsilon(t_nv_n)+\displaystyle\frac{1}{2}\int_{\mathbb{R}^3}(V^\beta_\varepsilon(x)-V(\varepsilon x)
(t_nv_n)^2dx\\
&+\displaystyle\frac{1}{p}\int_{\mathbb{R}^3}(P(\varepsilon x)-P^d_\varepsilon(x))
|t_nv_n|^pdx.
\end{aligned}
\end{equation*}
For
\begin{equation*}
\int_{\mathbb{R}^3}(V^\beta_\varepsilon(x)-V(\varepsilon x))
(t_nv_n)^2dx=\int_{\{x:V(\varepsilon x)<\beta\}}(\beta-V(\varepsilon x))
(t_nv_n)^2dx.
\end{equation*}
By $(PQ2)$, we know ${\{x:V(\varepsilon x)<\beta\}}$ is bounded. If $v_\varepsilon=0$, we have $v_n\rightarrow 0$ in $L^2_{loc}(\mathbb{R}^3)$. Then, by Lebesgue Dominated Convergence Theorem, we have
\begin{equation*}
\int_{\mathbb{R}^3}(V^\beta_\varepsilon(x)-V(\varepsilon x))
f^2(t_nv_n)dx=o_n(1).
\end{equation*}
For
\begin{equation*}
\int_{\mathbb{R}^3}(P(\varepsilon x)-P^b_\varepsilon(x))
|t_nv_n|^pdx=\int_{\{x:P(\varepsilon x)>b\}}(P(\varepsilon x)-b)
|t_nv_n|^pdx.
\end{equation*}
By the choice of $b$, we know ${\{x:P(\varepsilon x)>b\}}$ is bounded. Similarly, we also have
\begin{equation*}
\int_{\mathbb{R}^3}(P(\varepsilon x)-P^b_\varepsilon(x))
|t_nv_n|^pdx=o_n(1).
\end{equation*}
Thus we have \textbf{Claim 3}.
\par
Therefore,
\begin{equation*}
c^{\beta dQ_{max}}_\varepsilon\leq J^{\beta dQ_{max}}_\varepsilon(t_nv_n)=J_\varepsilon(t_nv_n)+o_n(1)\leq J_\varepsilon(v_n)+o_n(1).
\end{equation*}
Letting $n\rightarrow\infty$, we have $c^{\beta dQ_{max}}_\varepsilon\leq c_\varepsilon.$  By Lemma \ref{cv(y)}, we can obtain
\begin{equation*}
\limsup_{\varepsilon\to 0}c_\varepsilon\leq m^*_{V(0) P(0)Q(0)}=m^*_{\beta P_{max}Q_{max}}.
\end{equation*}
Then we have $\displaystyle\limsup_{\varepsilon\to 0}c^{\beta dQ_{max}}_\varepsilon\leq m^*_{\beta P_{max}Q_{max}}$. It follows from Lemma \ref{cbigm} that
\begin{equation*}
\limsup_{\varepsilon\to 0}c^{\beta dQ_{max}}_\varepsilon=\limsup_{\varepsilon\to 0}c^{V(0) dQ_{max}}_\varepsilon\geq m^*_{V(0) dQ_{max}}=m^*_{\beta dQ_{max}}.
\end{equation*}
Therefore, $m^*_{\beta dQ_{max}}\leq m^*_{\beta P_{max}Q_{max}}$, which contradicts Lemma \ref{le:comparel}. Thus, we have $v_\varepsilon\neq 0.$
\par
Moreover,
\begin{equation*}
\begin{aligned}
c_\varepsilon&=\displaystyle\lim_{n\to\infty}\Big[J_\varepsilon(v_n)-\displaystyle\frac{1}{4}\langle J_\varepsilon'(v_n),v_n\rangle\Big]\\
&=\displaystyle\lim_{n\to\infty}\Big[\frac{1}{4}\|v_n\|^2_\varepsilon+(\frac{1}{4}-\frac{1}{p})\int_{\mathbb{R}^3}P(\varepsilon x)|v_n|^pdx+(\frac{1}{4}-\frac{1}{6})\int_{\mathbb{R}^3}Q(\varepsilon x)|v_n|^6dx\Big]\\
&\geq\displaystyle\frac{1}{4}\|v_\varepsilon\|^2_\varepsilon+(\frac{1}{4}-\frac{1}{p})\int_{\mathbb{R}^3}P(\varepsilon x)|v_\varepsilon|^pdx+(\frac{1}{4}-\frac{1}{6})\int_{\mathbb{R}^3}Q(\varepsilon x)|v_\varepsilon|^6dx\\
&=J_\varepsilon(v_\varepsilon)-\displaystyle\frac{1}{4}\langle J_\varepsilon'(v_\varepsilon),v_\varepsilon\rangle\\
&=J_\varepsilon(v_\varepsilon).
\end{aligned}
\end{equation*}
Hence, $J_\varepsilon(v_\varepsilon)=c_\varepsilon.$ A similar argument to the one used in Lemma \ref{le:constant-achieved} shows $v_\varepsilon\in C^2(\mathbb{R}^3)$ and $v_\varepsilon>0$.
Then, $v_\varepsilon$ is a positive ground state solution of \eqref{transeq1}. Thus
$u_\varepsilon(x):=v_\varepsilon(\frac{x}{\varepsilon})$ a positive ground state solution of \eqref{maineq1}.
\hfill{$\Box$}

\begin{lemma}\label{lemma4.2}  Assume additionally that the potential functions $V(x)$, $P(x)$ and $Q(x)$ are uniformly continuous on $\mathbb{R}^3$.  Let $v_n:=v_{\varepsilon_n}$ be the solution obtained in Lemma \ref{lemma4.1} with $\varepsilon_n\rightarrow0$, as $n\to +\infty$. Then
\begin{itemize}
\item[$(1)$] there exists $y_n \in \mathbb{R}^3$ satisfying $\displaystyle\lim_{n\rightarrow\infty}dist(\varepsilon_ny_n,\mathcal{A}_V)=0$.
\item[$(2)$] up to a subsequence, $\displaystyle\lim_{n\rightarrow\infty}\varepsilon_ny_n=y_0.$  Set $\widetilde{v} _n(x)=v_n(x+y_n)$, then $\widetilde{v}_n(x)\rightarrow v$ in $E$, where $v$ is a positive ground state solution of
\begin{equation*}
-(a+b\int_{\mathbb{R}^3}|\nabla v|^2dx)\Delta v+V(y_0)v=P(y_0)|v|^{p-2}v+Q(y_0)|v|^4v,\quad x\in\mathbb{R}^3.
\end{equation*}
\end{itemize}
\end{lemma}
{\bf Proof:} Let $v_n$ be the positive ground state solution obtained in Lemma \ref{lemma4.1} with $\varepsilon_n\rightarrow0$. We claim that there exist $(y_n)\subset \mathbb{R}^3$ and $\rho,\eta>0$ such that
\begin{equation}\label{lionslemma1}
\displaystyle\limsup_{n\to\infty}\int_{B_\rho(y_n)}|v_n|^2dx\geq \eta.
\end{equation}
Suppose by contradiction that \eqref{lionslemma1} does not hold. Then, by Lions's Lemma, we can obtain
\begin{equation*}
v_n\rightarrow 0 \ \hbox{in}\  L^q(\mathbb{R}^3),\ \hbox{for}\ 2<q<6.
\end{equation*}
From
\begin{equation*}
\begin{aligned}
c_{\varepsilon_n}&=J_{\varepsilon_n}(v_n)-\displaystyle \frac{1}{p}\langle J'_{\varepsilon_n}(v_n),v_n\rangle\\
&=\displaystyle(\frac{1}{2}-\frac{1}{p})\|v_n\|^2_{\varepsilon_n}+(\frac{1}{2}-\frac{1}{p})b(\int_{\mathbb{R}^3}|\nabla v_n|^2dx)^2+(\frac{1}{p}-\frac{1}{6})\int_{\mathbb{R}^3}Q({\varepsilon_n} x)|v_n|^6dx,
\end{aligned}
\end{equation*}
we have
\begin{equation*}
(\frac{1}{2}-\frac{1}{p})\|v_n\|^2_{\varepsilon_n}\leq c_{\varepsilon_n}.
\end{equation*}
Then, by Lemma \ref{cv(y)} and $V(x)$ has a positive lower bound, it is not difficult to prove that $\{v_n\}$ is bounded in $E.$
\par
It follows from  $J'_{\varepsilon_n}(v_n)=0$ that
\begin{equation*}
\begin{aligned}
0&=\displaystyle\langle J_{\varepsilon_n}'(v_n),v_n\rangle\\
&=\displaystyle\|v_n\|^2_{\varepsilon_n}+b(\int_{\mathbb{R}^3}|\nabla v_n|^2dx)^2-\int_{\mathbb{R}^3}P({\varepsilon_n} x)|v_n|^{p}dx-\int_{\mathbb{R}^3}Q(\varepsilon_n x)|v_n|^{6}dx.
\end{aligned}
\end{equation*}
Thus,
\begin{equation*}
o(1)=\|v_n\|^2_{\varepsilon_n}+b(\int_{\mathbb{R}^3}|\nabla v_n|^2dx)^2-\int_{\mathbb{R}^3}Q(\varepsilon_n x)|v_n|^{6}dx.
\end{equation*}
Since $\{v_n\}$ is bounded in $E$, as $n\rightarrow\infty$, we can assume
\begin{equation}\label{lim13}
\|v_n\|^2_{\varepsilon_n}\rightarrow l_1,
\end{equation}
\begin{equation}\label{lim23}
b(\int_{\mathbb{R}^3}|\nabla v_n|^2dx)^2\rightarrow l_2,
\end{equation}
and
\begin{equation}\label{lim33}
\int_{\mathbb{R}^3}Q(\varepsilon_n x)|v_n|^{6}dx\rightarrow l_3.
\end{equation}
Then we have $l_3=l_1+l_2$.
\par
If $l_3=0$, we have $v_n\rightarrow 0$ in $E.$ Then $J_\varepsilon(v_n)\rightarrow 0,$ which contradicts Lemma \ref{le:lowboundedinN}. Thus, $l_3\neq0$.
\par
By the definition of the best constant $S$, we have
\begin{equation*}
S\leq \frac{\displaystyle\frac{1}{a}\int_{\mathbb{R}^3}a|\nabla v_n|^2dx}{\displaystyle(Q^{-1}_{max}\int_{\mathbb{R}^3}Q(\varepsilon x)|v_n|^{6}dx)^{1/3}}
\leq \frac{\displaystyle\frac{1}{a}\|v_n\|^2_{\varepsilon_n}}{\displaystyle(Q^{-1}_{max}\int_{\mathbb{R}^3}Q(\varepsilon x)|v_n|^{6}dx)^{1/3}}
\end{equation*}
and
\begin{equation*}
bS^2\leq \frac{ \displaystyle b(\int_{\mathbb{R}^3}|\nabla v_n|^2dx)^2}{\displaystyle Q^{-2/3}_{max}(\int_{\mathbb{R}^3}Q(\varepsilon x)|v_n|^{6}dx)^{2/3}}
\end{equation*}
Letting $n\to\infty$ in the above two inequalities, we have
\begin{equation*}
aSl_3^{1/3}\leq Q_{max}^{1/3}l_1 \quad  \hbox{and} \quad  bS^2l_3^{2/3}\leq Q_{max}^{2/3}l_2.
\end{equation*}
Therefore,
\begin{equation*}
l_1\geq aSQ_{max}^{-\frac{1}{3}}(l_1+l_2)^{\frac{1}{3}}\quad  \hbox{and} \quad  l_2\geq bS^2Q_{max}^{-\frac{2}{3}}(l_1+l_2)^{\frac{2}{3}}.
\end{equation*}
By Lemma \ref{inequality}, we have
\begin{equation*}
\frac{1}{3}l_1+\frac{1}{12}l_2\geq c^*:=\frac{ab}{4Q_{max}}S^3+\frac{(b^2S^4+4Q_{max} aS)^{3/2}}{24Q_{max}^2}+\frac{b^3S^6}{24Q_{max}^2}.
\end{equation*}
On the other hand,
\begin{equation*}
\begin{aligned}
c_{\varepsilon_n}&=J_{\varepsilon_n}(v_n)\\
&=\displaystyle\frac{1}{2}\|v_n\|^2_{\varepsilon_n}+\frac{b}{4}(\int_{\mathbb{R}^3}|\nabla v_n|^2dx)^2-\frac{1}{p}\int_{\mathbb{R}^3}P(\varepsilon_n x)|v_n|^{p}dx-\frac{1}{6}\int_{\mathbb{R}^3}Q(\varepsilon_n x)|v_n|^{6}dx\\
&=\displaystyle\frac{1}{2}\|v_n\|^2_{\varepsilon_n}+\frac{b}{4}(\int_{\mathbb{R}^3}|\nabla v_n|^2dx)^2-\frac{1}{6}\int_{\mathbb{R}^3}Q(\varepsilon_n x)|v_n|^{6}dx+o(1)\\
&=\displaystyle\frac{1}{2}\|v_n\|^2_{\varepsilon_n}+\frac{b}{4}(\int_{\mathbb{R}^3}|\nabla v_n|^2dx)^2-\frac{1}{6}(\|v_n\|^2_{\varepsilon_n}+b(\int_{\mathbb{R}^3}|\nabla v_n|^2dx)^2)+o(1)\\
&=\displaystyle\frac{1}{3}\|v_n\|^2_{\varepsilon_n}+\frac{b}{12}(\int_{\mathbb{R}^3}|\nabla v_n|^2dx)^2+o(1).\\
\end{aligned}
\end{equation*}
Using \eqref{lim13} and \eqref{lim23} in the above expression, we obtain
\begin{equation*}
\lim_{n\to\infty}c_{\varepsilon_n}= \frac{1}{3}l_1+\frac{1}{12}l_2\geq c^*,
\end{equation*}
which contradicts Lemma \ref{le:c-level}. Therefore, \eqref{lionslemma1} holds.
\par
Set $\widetilde{v}_n(x)=v_n(x+y_n)$. Then  $\widetilde{v}_n$ is bounded in $E$. Thus there exists $\widetilde{v}\in E$ such that $\widetilde{v}_n\rightharpoonup \widetilde{v}$ in $E.$ By \eqref{lionslemma1}, we know $\widetilde{v}\neq 0.$
Let $\widetilde{V}_{\varepsilon_n}(x)=V(\varepsilon_n(x+y_n)),$ $\widetilde{P}_{\varepsilon_n}(x)=P(\varepsilon_n(x+y_n))$ and $\widetilde{Q}_{\varepsilon_n}(x)=Q(\varepsilon_n(x+y_n)).$ Then $\widetilde{v}_n$ solves the following problems separately
\begin{equation}\label{nequation}
-(a+b\int_{\mathbb{R}^3}|\nabla v|^2dx)\Delta v+\widetilde{V}_{\varepsilon_n}(x)v=\widetilde{P}_{\varepsilon_n}(x)|v|^{p-2}v+\widetilde{Q}_{\varepsilon_n}(x)|v|^4v,\quad x\in\mathbb{R}^3.
\end{equation}
The corresponding energy functional
\begin{equation*}
\begin{aligned}
\widetilde{J}_{\varepsilon_n}(v)=&\displaystyle\frac{1}{2}\int_{\mathbb{R}^3}(a|\nabla v|^2+\widetilde{V}_{\varepsilon_n}(x)v^2)dx+\frac{b}{4}(\int_{\mathbb{R}^3}|\nabla v|^2dx)^2\\
&\displaystyle-\frac{1}{p}\int_{\mathbb{R}^3}\widetilde{P}_{\varepsilon_n}(x)|v|^{p}dx-\frac{1}{6}\int_{\mathbb{R}^3}\widetilde{Q}_{\varepsilon_n}(x)|v|^{6}dx,
\end{aligned}
\end{equation*}
\par
\textbf{Claim 1:} $\varepsilon_ny_n$ must be bounded.
\par
Otherwise, without loss of generality, we assume $\varepsilon_ny_n\rightarrow\infty$ as $n\rightarrow\infty.$ Up to a subsequence, we have
\begin{equation*}
\begin{aligned}
&V(\varepsilon_ny_n)\rightarrow V_0\geq \beta,\\
&P(\varepsilon_ny_n)\rightarrow P_0< P_{max},\\
&Q(\varepsilon_ny_n)\rightarrow Q_0\leq Q_{max}.
\end{aligned}
\end{equation*}
Then $\widetilde{v}$ is a solution of the following equation
\begin{equation*}
-(a+b\int_{\mathbb{R}^3}|\nabla v|^2dx)\Delta v+V_0v=P_0|v|^{p-2}v+Q_0|v|^4v,\quad x\in\mathbb{R}^3.
\end{equation*}
In fact, for any test function $\phi\in C_0^\infty(\mathbb{R}^3)$, since $\widetilde{v}_n$ is a solution of equation \eqref{nequation}, we have
\begin{equation*}
\begin{aligned}
0&=\displaystyle\int_{\mathbb{R}^3}(a\nabla \widetilde{v}_n \nabla \phi+ \widetilde{V}_{\varepsilon_n}(x)\widetilde{v}_n \phi)dx+b\int_{\mathbb{R}^3}|\nabla \widetilde{v}_n|^2dx\int_{\mathbb{R}^3}\nabla \widetilde{v}_n \nabla \phi dx\\
&\quad\displaystyle-\int_{\mathbb{R}^3}\widetilde{P}_{\varepsilon_n}(x)|\widetilde{v}_n|^{p-2}\widetilde{v}_n\phi dx-\int_{\mathbb{R}^3}\widetilde{Q}_{\varepsilon_n}(x)|\widetilde{v}_n|^{4}\widetilde{v}_n\phi dx,
\end{aligned}
\end{equation*}
Assume that $\int_{\mathbb{R}^3}|\nabla \widetilde{v}_n|^2dx\to A^2$, as $n\to \infty.$ Note the uniformly continuity of $V(x),P(x),Q(x)$, we have
\begin{equation*}
\begin{aligned}
0&=\displaystyle\int_{\mathbb{R}^3}(a\nabla \widetilde{v} \nabla \phi+ V_0\widetilde{v} \phi)dx+bA^2\int_{\mathbb{R}^3}\nabla \widetilde{v} \nabla \phi dx\\
&\quad\displaystyle-\int_{\mathbb{R}^3} P_0|\widetilde{v}|^{p-2}\widetilde{v}\phi dx-\int_{\mathbb{R}^3}Q_0|\widetilde{v}|^{4}\widetilde{v}\phi dx.
\end{aligned}
\end{equation*}
Thus,we can get
\begin{equation*}
\int_{\mathbb{R}^3}(a|\nabla \widetilde{v}|^2+V_0\widetilde{v}^2)dx+ bA^2\int_{\mathbb{R}^3}|\nabla \widetilde{v}|^2dx
=\int_{\mathbb{R}^3} P_0|\widetilde{v}|^{p}dx+\int_{\mathbb{R}^3}Q_0|\widetilde{v}|^{6}dx.
\end{equation*}
It is easy to know that $\int_{\mathbb{R}^3}|\nabla \widetilde{v}|^2dx\leq A^2$. If $\int_{\mathbb{R}^3}|\nabla \widetilde{v}|^2dx< A^2$, we have
\begin{equation*}
\int_{\mathbb{R}^3}(a|\nabla \widetilde{v}|^2+V_0\widetilde{v}^2)dx+ b(\int_{\mathbb{R}^3}|\nabla \widetilde{v}|^2dx)^2
<\int_{\mathbb{R}^3} P_0|\widetilde{v}|^{p}dx+\int_{\mathbb{R}^3}Q_0|\widetilde{v}|^{6}dx.
\end{equation*}
Then, there exists $0<\widetilde{t}<1$ such that $\widetilde{t}\widetilde{v}\in N^*_{V_0P_0Q_0}$.  Thus,  we have
\begin{equation*}
\begin{aligned}
m^*_{V_0P_0Q_0}&=\Phi^{\ast}_{V_0P_0Q_0}(\widetilde{t}\widetilde{v})-\displaystyle\frac{1}{4}\langle (\Phi^{\ast}_{V_0P_0Q_0})'(\widetilde{t}\widetilde{v}),\widetilde{t}\widetilde{v}\rangle\\
&=\displaystyle\frac{\widetilde{t}^2}{4}\int_{\mathbb{R}^3}(a|\nabla \widetilde{v}|^2+V_0\widetilde{v}^2)dx+(\frac{\widetilde{t}^p}{4}-\frac{\widetilde{t}^p}{p})\int_{\mathbb{R}^3} P_0|\widetilde{v}|^{p}dx+(\frac{\widetilde{t}^6}{4}-\frac{\widetilde{t}^6}{6})\int_{\mathbb{R}^3}Q_0|\widetilde{v}|^{6}dx\\
&<\displaystyle\frac{1}{4}\int_{\mathbb{R}^3}(a|\nabla \widetilde{v}|^2+V_0\widetilde{v}^2)dx+(\frac{1}{4}-\frac{1}{p})\widetilde\int_{\mathbb{R}^3} P_0|\widetilde{v}|^{p}dx+(\frac{1}{4}-\frac{1}{6})\int_{\mathbb{R}^3}Q_0|\widetilde{v}|^{6}dx\\
&\leq\displaystyle\liminf_{n\to\infty}\Big[\frac{1}{4}\int_{\mathbb{R}^3}(a|\nabla \widetilde{v}_n|^2+\widetilde{V}_{\varepsilon_n}(x)\widetilde{v}_n^2)dx+(\frac{1}{4}-\frac{1}{p})\widetilde\int_{\mathbb{R}^3} \widetilde{P}_{\varepsilon_n}(x)|\widetilde{v}_n|^{p}dx\\
&\quad\displaystyle+(\frac{1}{4}-\frac{1}{6})\int_{\mathbb{R}^3}\widetilde{Q}_{\varepsilon_n}(x)|\widetilde{v}_n|^{6}dx\Big]\\
&=\displaystyle\liminf_{n\to\infty}\Big[\widetilde{J}_{\varepsilon_n}(\widetilde{v}_n)-\displaystyle\frac{1}{4}\langle \widetilde{J}_{\varepsilon_n}'(\widetilde{v}_n),\widetilde{v}_n\rangle\Big]
=\displaystyle\liminf_{n\to\infty}J_{\varepsilon_n}(v_n)
=\liminf_{n\to\infty}c_{\varepsilon_n}.
\end{aligned}
\end{equation*}
Then, by Lemma \ref{cv(y)},  we have
\begin{equation*}
m^*_{V_0P_0Q_0}\leq m^*_{V(0)P(0)Q(0)}=m^*_{\beta P_{max}Q_{max}},
\end{equation*}
which contradicts Lemma \ref{le:comparel}.
\par
Therefore,
\begin{equation*}
-(a+b\int_{\mathbb{R}^3}|\nabla \widetilde{v}|^2dx)\Delta \widetilde{v}+V_0\widetilde{v}=P_0|\widetilde{v}|^{p-2}\widetilde{v}+Q_0|\widetilde{v}|^4\widetilde{v},\quad x\in\mathbb{R}^3.
\end{equation*}
It follows that
\begin{equation*}
\begin{aligned}
m^*_{\beta P_{max}Q_{max}}&<m^*_{V_0P_0Q_0}\\
&\leq\Phi^{\ast}_{V_0P_0Q_0}(\widetilde{v})\\
&=\Phi^{\ast}_{V_0P_0Q_0}(\widetilde{v})-\displaystyle\frac{1}{4}\langle (\Phi^{\ast}_{V_0P_0Q_0})'(\widetilde{v}),\widetilde{v}\rangle\\
&=\displaystyle\frac{1}{4}\int_{\mathbb{R}^3}(a|\nabla \widetilde{v}|^2+V_0\widetilde{v}^2)dx+(\frac{1}{4}-\frac{1}{p})\int_{\mathbb{R}^3} P_0|\widetilde{v}|^{p}dx+(\frac{1}{4}-\frac{1}{6})\int_{\mathbb{R}^3}Q_0|\widetilde{v}|^{6}dx\\
&\leq\displaystyle\liminf_{n\to\infty}\Big[\frac{1}{4}\int_{\mathbb{R}^3}(a|\nabla \widetilde{v}_n|^2+\widetilde{V}_{\varepsilon_n}(x)\widetilde{v}_n^2)dx+(\frac{1}{4}-\frac{1}{p})\int_{\mathbb{R}^3} \widetilde{P}_{\varepsilon_n}(x)|\widetilde{v}_n|^{p}dx\\
&\quad\displaystyle+(\frac{1}{4}-\frac{1}{6})\int_{\mathbb{R}^3}\widetilde{Q}_{\varepsilon_n}(x)|\widetilde{v}_n|^{6}dx\Big]\\
&=\displaystyle\liminf_{n\to\infty}\Big[\widetilde{J}_{\varepsilon_n}(\widetilde{v}_n)-\displaystyle\frac{1}{4}\langle \widetilde{J}_{\varepsilon_n}'(\widetilde{v}_n),\widetilde{v}_n\rangle\Big]
=\displaystyle\liminf_{n\to\infty}J_{\varepsilon_n}(v_n)
=\liminf_{n\to\infty}c_{\varepsilon_n},
\end{aligned}
\end{equation*}
which contradicts Lemma \ref{cv(y)}. Thus $\varepsilon_ny_n$ must be bounded. And, up to a subsequence, we can assume $\varepsilon_ny_n\rightarrow y_0$.
\par
\textbf{Claim 2:} $y_0\in \mathcal{A}_V.$
\par
If $y_0\notin \mathcal{A}_V,$ we have two cases.
\begin{itemize}
\item[$(1)$] $\beta<V(y_0)$, $P(y_0)=P_{max}$ and $Q(y_0)=Q_{max}$, then $m^*_{\beta P_{max}Q_{max}}<m^*_{V(y_0)P(y_0)Q(y_0)}$.
\item[$(2)$] $\beta\leq V(y_0)$, $P(y_0)<P_{max}$ or $Q(y_0)<Q_{max}$, then  $m^*_{\beta P_{max}Q_{max}}<m^*_{V(y_0)P(y_0)Q(y_0)}$.
\end{itemize}
From \textbf{Claim 1}, we know that $\widetilde{v}$ is a solution of the following equation
\begin{equation}\label{nequation2}
-(a+b\int_{\mathbb{R}^3}|\nabla v|^2dx)\Delta v+V(y_0)v=P(y_0)|v|^{p-2}v+Q(y_0)|v|^4v,\quad x\in\mathbb{R}^3.
\end{equation}
Similar to the arguments in \textbf{Claim 1}, we have $m^*_{V(y_0)P(y_0)Q(y_0)}\leq\displaystyle\liminf_{n\to\infty}c_{\varepsilon_n}.$
\par
Applying Lemma \ref{cv(y)} and Lemma \ref{le:comparel}, we get
\begin{equation*}
\limsup_{n\rightarrow \infty}c_{\varepsilon_n}\leq m^*_{\beta P_{max}Q_{max}}<m^*_{V(y_0)P(y_0)Q(y_0)}\leq\liminf_{n\to\infty}c_{\varepsilon_n},
\end{equation*}
which is absurd. Therefore, $y_0\in \mathcal{A}_V,$ which means that
\begin{equation*}
\lim_{n\rightarrow\infty}dist(\varepsilon_ny_n,\mathcal{A}_V)=0.
\end{equation*}
\par
\textbf{Claim 3:} $\widetilde{v}$ is a positive ground state solution of \eqref{nequation2}.
\par
Repeating the arguments in \textbf{Claim 1} again, we get
\begin{equation*}
\Phi^{\ast}_{V(y_0)P(y_0)Q(y_0)}(\widetilde{v})\leq\displaystyle\liminf_{n\to\infty}J_{\varepsilon_n}(v_n)
\leq\liminf_{n\to\infty}c_{\varepsilon_n}\leq\limsup_{n\to\infty}c_{\varepsilon_n}\leq m^*_{V(y_0)P(y_0)Q(y_0)}.
\end{equation*}
So we get
\begin{equation*}
\Phi^{\ast}_{V(y_0)P(y_0)Q(y_0)}(\widetilde{v})=m^*_{V(y_0)P(y_0)Q(y_0)}.
\end{equation*}
Thus $\widetilde{v}$ is a ground state solution. By the theory of
elliptic regularity, $\widetilde{v}\in C^2(\mathbb{R}^3)$, and by using strong maximum principle, we get $\widetilde{v}>0$ in $\mathbb{R}^3$.
\par
\textbf{Claim 4:} $\widetilde{v}_n$ converges strongly to $\widetilde{v}$ in $E$.
\par
From \textbf{Claim 3}, we know that
\begin{equation*}
\lim_{n\to\infty}\widetilde{J}_{\varepsilon_n}(\widetilde{v}_n)=\lim_{n\to\infty}J_{\varepsilon_n}(v_n)
=\lim_{n\to\infty}c_{\varepsilon_n}=\Phi^{\ast}_{V(y_0)P(y_0)Q(y_0)}(\widetilde{v}).
\end{equation*}
It follows that
\begin{equation*}
\begin{aligned}
m^*_{V(y_0)P(y_0)Q(y_0)}&=\Phi^{\ast}_{V(y_0)P(y_0)Q(y_0)}(\widetilde{v})-\displaystyle\frac{1}{4}\langle (\Phi^{\ast}_{V(y_0)P(y_0)Q(y_0)})'(\widetilde{v}),\widetilde{v}\rangle\\
&=\displaystyle\frac{1}{4}\int_{\mathbb{R}^3}(a|\nabla \widetilde{v}|^2+V(y_0)\widetilde{v}^2)dx+(\frac{1}{4}-\frac{1}{p})\int_{\mathbb{R}^3} P(y_0)|\widetilde{v}|^{p}dx\\
&\quad\displaystyle+(\frac{1}{4}-\frac{1}{6})\int_{\mathbb{R}^3}Q(y_0)|\widetilde{v}|^{6}dx\\
&\leq\displaystyle\liminf_{n\to\infty}\Big[\frac{1}{4}\int_{\mathbb{R}^3}(a|\nabla \widetilde{v}_n|^2+\widetilde{V}_{\varepsilon_n}(x)\widetilde{v}_n^2)dx+(\frac{1}{4}-\frac{1}{p})\int_{\mathbb{R}^3} \widetilde{P}_{\varepsilon_n}(x)|\widetilde{v}_n|^{p}dx\\
&\quad\displaystyle+(\frac{1}{4}-\frac{1}{6})\int_{\mathbb{R}^3}\widetilde{Q}_{\varepsilon_n}(x)|\widetilde{v}_n|^{6}dx\Big]\\
&=\displaystyle\liminf_{n\to\infty}\Big[\widetilde{J}_{\varepsilon_n}(\widetilde{v}_n)-\displaystyle\frac{1}{4}\langle \widetilde{J}_{\varepsilon_n}'(\widetilde{v}_n),\widetilde{v}_n\rangle\Big]\\
&=\Phi^{\ast}_{V(y_0)P(y_0)Q(y_0)}(\widetilde{v}).
\end{aligned}
\end{equation*}
Thus, as $n\to \infty$, we have
\begin{equation}\label{nlim}
\int_{\mathbb{R}^3}(a|\nabla \widetilde{v}_n|^2+\widetilde{V}_{\varepsilon_n}(x)\widetilde{v}_n^2)dx\to \int_{\mathbb{R}^3}(a|\nabla \widetilde{v}|^2+V(y_0)\widetilde{v}^2)dx.
\end{equation}
Since
\begin{equation*}
\int_{\mathbb{R}^3}a|\nabla \widetilde{v}|^2dx\leq\liminf_{n\to\infty}\int_{\mathbb{R}^3}a|\nabla \widetilde{v}_n|^2dx
\end{equation*}
and
\begin{equation*}
\int_{\mathbb{R}^3}V(y_0)\widetilde{v}^2dx\leq \liminf_{n\to\infty}\int_{\mathbb{R}^3}\widetilde{V}_{\varepsilon_n}(x)\widetilde{v}_n^2)dx,
\end{equation*}
it follows from \eqref{nlim} that
\begin{equation*}
\lim_{n\to\infty}\int_{\mathbb{R}^3}a|\nabla \widetilde{v}_n|^2dx=\int_{\mathbb{R}^3}a|\nabla \widetilde{v}|^2dx
\end{equation*}
and
\begin{equation*}
\lim_{n\to\infty}\int_{\mathbb{R}^3}\widetilde{V}_{\varepsilon_n}(x)\widetilde{v}_n^2dx=\int_{\mathbb{R}^3}V(y_0)\widetilde{v}^2dx.
\end{equation*}
Then, it is easy to prove that
\begin{equation*}
\int_{\mathbb{R}^3}(a|\nabla \widetilde{v}_n|^2+V(y_0)\widetilde{v}_n^2)dx\to \int_{\mathbb{R}^3}(a|\nabla \widetilde{v}|^2+V(y_0)\widetilde{v}^2)dx.
\end{equation*}
Thus, we have $\|\widetilde{v}_n\|\to\|\widetilde{v}\|$, noting that $\widetilde{v}_n\rightharpoonup \widetilde{v}$ in $E$, So  $\widetilde{v}_n\rightarrow\widetilde{v}$ in $E$ is obtained.
\hfill{$\Box$}

\begin{remark}\label{re1}
In fact, from the proof of Lemma \ref{lemma4.2}, we can get the following results.
\begin{itemize}
\item[$(1)$] There exists $\varepsilon^*>0$, a family $\{y_\varepsilon\}\subset \mathbb{R}^3$ and $\rho,\eta>0$ such that, for all $\varepsilon\in(0,\varepsilon^*)$,
\begin{equation}\label{lionslemma12}
\begin{array}{ll}
\displaystyle\int_{B_\rho(y_\varepsilon)}|v_\varepsilon|^2dx\geq \eta.
\end{array}
\end{equation}
\item[$(2)$] $\{\varepsilon y_\varepsilon\}$ is bounded, satisfying $\displaystyle\lim_{\varepsilon\rightarrow 0}dist(\varepsilon y_\varepsilon,\mathcal{A}_V)=0$.
\end{itemize}
\end{remark}

\begin{lemma}\label{lemma4.3}  There exists $\varepsilon^*>0$ such that
\begin{equation*}
\displaystyle\lim_{|x|\rightarrow\infty}\widetilde{v}_\varepsilon(x)=0\quad \text{uniformly on}\quad \varepsilon\in (0,\varepsilon^*),
\end{equation*}
and there exists $C>0$ independent of $\varepsilon$ such that $|\widetilde{v}_\varepsilon|_\infty\leq C$ uniformly on $\varepsilon\in (0,\varepsilon^*)$, where $\widetilde{v}_\varepsilon$ are obtained in Lemma \ref{lemma4.2}. Furthermore, there exist constants $C,c>0$ such that
\begin{equation*}
|\widetilde{v}_\varepsilon(x)|\leq Cexp(-c|x|)
\end{equation*}
for all $x\in \mathbb{R}^3$.
\end{lemma}
{\bf Proof:} The proof of this lemma can be obtained from Lemma 4.4 and Lemma 4.5 in \cite{MR3347410}.
\hfill{$\Box$}

\begin{remark}\label{re2}
From \eqref{lionslemma12} and Lemma \ref{lemma4.3}, we have
\begin{equation*}
\frac{\eta}{2}\leq\int_{B_\rho(0)}|\widetilde{v}_\varepsilon|^2dx\leq C|\widetilde{v}_\varepsilon|_\infty.
\end{equation*}
Thus, there exists $\eta'>0$, such that $|\widetilde{v}_\varepsilon|_\infty\geq \eta'.$  If $b_\varepsilon$ is a maximum point of $\widetilde{v}_\varepsilon$, by $\displaystyle\lim_{|x|\rightarrow\infty}\widetilde{v}_\varepsilon(x)=0$ uniformly on $\varepsilon\in (0,\varepsilon^*)$, we can get $R_0>0$ such that $|b_\varepsilon|\leq R_0.$
\end{remark}
\par
{\bf The proof of Theorem 1.1:}
\par
By Lemma \ref{lemma4.1}, for $\varepsilon>0$ small enough, problem \eqref{maineq1} has a positive ground state solution  $u_\varepsilon(x)=v_\varepsilon(\frac{x}{\varepsilon})$.
\par
$(1)$ From Remark \ref{re2}, $\widetilde{v}_\varepsilon$ has a maximum point $b_\varepsilon$. Then $v_\varepsilon$ has a maximum point $z_\varepsilon:=b_\varepsilon+y_\varepsilon$. Thus, $u_\varepsilon(x)$ has maximum value at $x_\varepsilon:=\varepsilon z_\varepsilon$. Noting the boundness of $b_\varepsilon$, by Remark \ref{re1}, we have $\displaystyle\lim_{\varepsilon\rightarrow0}dist(x_\varepsilon,\mathcal{A}_V)=0.$ Moreover,
it follows from Lemma \ref{lemma4.3} that
\begin{equation*}
u_\varepsilon(x)=v_\varepsilon(\frac{x}{\varepsilon})=\widetilde{v}_\varepsilon(\frac{x}{\varepsilon}-y_\varepsilon)\leq Cexp(-c |\frac{x}{\varepsilon}-y_\varepsilon|)\leq Cexp(-\frac{c}{\varepsilon} |x-x_\varepsilon|).
\end{equation*}
\par
$(2)$ Since $\widetilde{x}_\varepsilon$ is a maximum point of $u_\varepsilon,$ then $\widetilde{b}_\varepsilon:=\frac{\widetilde{x}_\varepsilon}{\varepsilon}-y_\varepsilon$ is the maximum point of $\widetilde{v}_\varepsilon.$ In view of Remark \ref{re2}, we know $\widetilde{b}_\varepsilon$ is bounded. Moreover,  $\varepsilon(\widetilde{b}_\varepsilon+y_\varepsilon)=\widetilde{x}_\varepsilon\rightarrow x_0$ as $\varepsilon\rightarrow0.$
\par
On the other hand, by \eqref{lionslemma12}, there exist $\rho,\eta>0$ such that
\begin{equation*}
\displaystyle\limsup_{\varepsilon\to0}\int_{B_\rho(y_\varepsilon)}|v_\varepsilon|^2dx\geq \eta.
\end{equation*}
So we have
\begin{equation*}
\displaystyle\limsup_{\varepsilon\to0}\int_{B_{\rho+R_0}(y_\varepsilon+\widetilde{b}_\varepsilon)}|v_\varepsilon|^2dx\geq \displaystyle\limsup_{\varepsilon\to0}\int_{B_\rho(y_\varepsilon)}|v_\varepsilon|^2dx\geq \eta.
\end{equation*}
Then, using the same argument as in the proof of Lemma \ref{lemma4.2}, we get $v_\varepsilon(x+\widetilde{b}_\varepsilon+y_\varepsilon)\rightarrow v$ in $E$, as $\varepsilon\rightarrow0,$ where $v$ is a positive ground state solution of
\begin{equation*}
-(a+b\int_{\mathbb{R}^3}|\nabla v|^2dx)\Delta v+V(y_0)v=P(y_0)|v|^{p-2}v+Q(y_0)|v|^4v,\quad x\in\mathbb{R}^3.
\end{equation*}
Thus,
$u_\varepsilon(\varepsilon x+\widetilde{x}_\varepsilon)=v_\varepsilon(x+\widetilde{b}_\varepsilon+y_\varepsilon)\rightarrow v$ in $E$, as $\varepsilon\rightarrow 0$.
\hfill{$\Box$}
\par
The proof of Theorem 1.2 is similar to Theorem 1.1, so we omit the detail.


\section*{Acknowledgments}
We would like to thank the anonymous referee for his/her careful readings of our manuscript and the useful comments made for its improvement. The first author thanks his advisor Prof. Zhongwei Tang for suggestions and help.  The second author also thanks the support of RTG 2419 by the German Science Foundation (DFG).

\bibliographystyle{plain}
\bibliography{Kirchhoff}

\begin{thebibliography}{10}

\bibitem{Alves-Correa-Ma2005CMP}
C.~O. Alves, F.~J. S.~A. Corr\^ea, and T.~F. Ma.
\newblock Positive solutions for a quasilinear elliptic equation of {K}irchhoff
  type.
\newblock {\em Comput. Math. Appl.}, 49(1):85--93, 2005.

\bibitem{Arosio-Panizzi1996TAMS}
A.~Arosio and S.~Panizzi.
\newblock On the well-posedness of the {K}irchhoff string.
\newblock {\em Trans. Amer. Math. Soc.}, 348(1):305--330, 1996.

\bibitem{Bernstein1940BASUS}
S.~Bernstein.
\newblock Sur une classe d'\'equations fonctionnelles aux d\'eriv\'ees
  partielles.
\newblock {\em Bull. Acad. Sci. URSS. S\'er. Math. [Izvestia Akad. Nauk SSSR]},
  4:17--26, 1940.

\bibitem{Chen-Kuo-Wu2011JDE}
C.~Chen, Y.~Kuo, and T.~Wu.
\newblock The {N}ehari manifold for a {K}irchhoff type problem involving
  sign-changing weight functions.
\newblock {\em J. Differential Equations}, 250(4):1876--1908, 2011.

\bibitem{DpAncona-Spagnolo1992IM}
P.~D'Ancona and S.~Spagnolo.
\newblock Global solvability for the degenerate {K}irchhoff equation with real
  analytic data.
\newblock {\em Invent. Math.}, 108(2):247--262, 1992.

\bibitem{Deng-Peng-Shuai2015JFA}
Y.~Deng, S.~Peng, and W.~Shuai.
\newblock Existence and asymptotic behavior of nodal solutions for the
  {K}irchhoff-type problems in {$\Bbb{R}^3$}.
\newblock {\em J. Funct. Anal.}, 269(11):3500--3527, 2015.

\bibitem{Ding-Liu2012JDE}
Y.~Ding and X.~Liu.
\newblock Semi-classical limits of ground states of a nonlinear {D}irac
  equation.
\newblock {\em J. Differential Equations}, 252(9):4962--4987, 2012.

\bibitem{Ding-Liu2013MM}
Y.~Ding and X.~Liu.
\newblock Semiclassical solutions of {S}chr\"{o}dinger equations with magnetic
  fields and critical nonlinearities.
\newblock {\em Manuscripta Math.}, 140(1-2):51--82, 2013.

\bibitem{He-Zou2009NA}
X.~He and W.~Zou.
\newblock Infinitely many positive solutions for {K}irchhoff-type problems.
\newblock {\em Nonlinear Anal.}, 70(3):1407--1414, 2009.

\bibitem{He-Zou2012JDE}
X.~He and W.~Zou.
\newblock Existence and concentration behavior of positive solutions for a
  {K}irchhoff equation in {$\Bbb R^3$}.
\newblock {\em J. Differential Equations}, 252(2):1813--1834, 2012.

\bibitem{He-Zou2014ADM}
X.~He and W.~Zou.
\newblock Ground states for nonlinear {K}irchhoff equations with critical
  growth.
\newblock {\em Ann. Mat. Pura Appl. (4)}, 193(2):473--500, 2014.

\bibitem{He-Li2015CVPDE}
Y.~He and G.~Li.
\newblock Standing waves for a class of {K}irchhoff type problems in
  {$\Bbb{R}^3$} involving critical {S}obolev exponents.
\newblock {\em Calc. Var. Partial Differential Equations}, 54(3):3067--3106,
  2015.

\bibitem{He-Li-Peng2014ANS}
Y.~He, G.~Li, and S.~Peng.
\newblock Concentrating bound states for {K}irchhoff type problems in {$\Bbb
  R^3$} involving critical {S}obolev exponents.
\newblock {\em Adv. Nonlinear Studies}, 14(2):483--510, 2014.

\bibitem{Kirchhoff1883}
G.~Kirchhoff.
\newblock {\em Vorlesungen \"uber Mechanik}.
\newblock Birkh\"auser Basel, 1883.

\bibitem{MR4021897}
G.~Li, P.~Luo, S.~Peng, C.~Wang, and C.~Xiang.
\newblock A singularly perturbed {K}irchhoff problem revisited.
\newblock {\em J. Differential Equations}, 268(2):541--589, 2020.

\bibitem{Li-Ye2014JDE}
G.~Li and H.~Ye.
\newblock Existence of positive ground state solutions for the nonlinear
  {K}irchhoff type equations in {$\Bbb R^3$}.
\newblock {\em J. Differential Equations}, 257(2):566--600, 2014.

\bibitem{Li-YeMMAS}
G.~Li and H.~Ye.
\newblock Existence of positive solutions for nonlinear {K}irchhoff type
  problems in {$\Bbb R^3$} with critical {S}obolev exponent.
\newblock {\em Math. Methods Appl. Sci.}, 37(16):2570--2584, 2014.

\bibitem{Lions1978NHMS}
J.~L. Lions.
\newblock On some questions in boundary value problems of mathematical physics.
\newblock 30:284--346, 1978.

\bibitem{MR3347410}
Z.~Liu and S.~Guo.
\newblock Existence and concentration of positive ground states for a
  {K}irchhoff equation involving critical {S}obolev exponent.
\newblock {\em Z. Angew. Math. Phys.}, 66(3):747--769, 2015.

\bibitem{Ma-Rivera2003AML}
T.~F. Ma and J.~E. Mu\~noz Rivera.
\newblock Positive solutions for a nonlinear nonlocal elliptic transmission
  problem.
\newblock {\em Appl. Math. Lett.}, 16(2):243--248, 2003.

\bibitem{Mao-Zhang2009NL}
A.~Mao and Z.~Zhang.
\newblock Sign-changing and multiple solutions of {K}irchhoff type problems
  without the {P}.{S}. condition.
\newblock {\em Nonlinear Anal.}, 70(3):1275--1287, 2009.

\bibitem{Ono1997JDE}
K.~Ono.
\newblock Global existence, decay, and blowup of solutions for some mildly
  degenerate nonlinear {K}irchhoff strings.
\newblock {\em J. Differential Equations}, 137(2):273--301, 1997.

\bibitem{Perera-Zhang2006JDE}
K.~Perera and Z.~Zhang.
\newblock Nontrivial solutions of {K}irchhoff-type problems via the {Y}ang
  index.
\newblock {\em J. Differential Equations}, 221(1):246--255, 2006.

\bibitem{Pohozaev1975MS}
S.~I. Pohozaev.
\newblock A certain class of quasilinear hyperbolic equations.
\newblock {\em Mat. Sb. (N.S.)}, 96(138):152--166, 168, 1975.

\bibitem{Shuai2015JDE}
W.~Shuai.
\newblock Sign-changing solutions for a class of {K}irchhoff-type problem in
  bounded domains.
\newblock {\em J. Differential Equations}, 259(4):1256--1274, 2015.

\bibitem{Wang-Tian-Xu-Zhang2012JDE}
J.~Wang, L.~Tian, J.~Xu, and F.~Zhang.
\newblock Multiplicity and concentration of positive solutions for a
  {K}irchhoff type problem with critical growth.
\newblock {\em J. Differential Equations}, 253(7):2314--2351, 2012.

\bibitem{Wang1993CMP}
X.~Wang.
\newblock On concentration of positive bound states of nonlinear
  {S}chr\"{o}dinger equations.
\newblock {\em Comm. Math. Phys.}, 153(2):229--244, 1993.

\bibitem{Wang-Zeng1997SIAMJMA}
X.~Wang and B.~Zeng.
\newblock On concentration of positive bound states of nonlinear
  {S}chr\"{o}dinger equations with competing potential functions.
\newblock {\em SIAM J. Math. Anal.}, 28(3):633--655, 1997.

\bibitem{Willem1996}
M.~Willem.
\newblock {\em Minimax theorems}, volume~24 of {\em Progress in Nonlinear
  Differential Equations and their Applications}.
\newblock Birkh\"auser Boston, Inc., Boston, MA, 1996.

\bibitem{Wu2011NA}
X.~Wu.
\newblock Existence of nontrivial solutions and high energy solutions for
  {S}chr\"odinger-{K}irchhoff-type equations in {${\bf R}^N$}.
\newblock {\em Nonlinear Anal. Real World Appl.}, 12(2):1278--1287, 2011.

\bibitem{MR3518335}
M.~Yang.
\newblock Concentration of positive ground state solutions for
  {S}chr\"{o}dinger-{M}axwell systems with critical growth.
\newblock {\em Adv. Nonlinear Stud.}, 16(3):389--408, 2016.

\bibitem{MR3994307}
Z.~Yang, Y.~Yu, and F.~Zhao.
\newblock The concentration behavior of ground state solutions for a critical
  fractional {S}chr\"{o}dinger-{P}oisson system.
\newblock {\em Math. Nachr.}, 292(8):1837--1868, 2019.

\bibitem{Zhang-Perera2006JMAA}
Z.~Zhang and K.~Perera.
\newblock Sign changing solutions of {K}irchhoff type problems via invariant
  sets of descent flow.
\newblock {\em J. Math. Anal. Appl.}, 317(2):456--463, 2006.

\end{thebibliography}

\end{document}